\documentclass[11pt]{article}
\usepackage[utf8]{inputenc}
\usepackage{graphicx}
\usepackage{color}
\usepackage{float}
\usepackage{todonotes}
\usepackage{amsfonts}
\usepackage{amsthm}
\usepackage{amsmath}
\usepackage{amssymb}
\usepackage{stmaryrd}
\usepackage[top=2.cm, bottom=2.cm, left=2.5cm, right=2.5cm]{geometry}
\usepackage[english]{babel}
\usepackage[T1]{fontenc}
\usepackage{hyperref}
\usepackage{slashbox}
\usepackage{algorithm}
\usepackage{algorithmicx}
\usepackage{algpseudocode}
\usepackage{float}
\usepackage{csquotes}
\usepackage{colortbl}
\usepackage{comment} 

\usepackage[backend=bibtex, style=alphabetic]{biblatex} 
\title{Numerical resolution of McKean-Vlasov FBSDEs using neural networks 
\thanks{This work is supported by  FiME, Laboratoire de Finance des March\'es de l'Energie.}}

\author{Maximilien \textsc{Germain}
\footnote{EDF R\&D, Université de Paris, LPSM  \sf \href{mailto:Maximilien.Germain at edf.fr}{mgermain at lpsm.fr}}, Joseph \textsc{Mikael}
\footnote{EDF R\&D  \sf \href{mailto:Joseph.Mikael at edf.fr}{joseph.mikael at edf.fr}}, Xavier \textsc{Warin}
\footnote{EDF R\&D \& FiME \sf \href{mailto:xavier.warin at edf.fr}{xavier.warin at edf.fr}} }

\date{\today\\ {\it to appear in Methodology and Computing in Applied Probability}}

\newcommand{\norme}[1]{\left\Vert #1\right\Vert}
\newcommand{\Prob}{\mathbb{P}}
\newcommand{\R}{\mathbb{R}}

\newcommand{\E}{\mathbb{E}}

\newcommand{\di}{\mathrm{d}}
\newtheorem{Rem}{Remark}
\def\Pc{{\cal P}}
\def\Lc{{\cal L}}

\DeclareMathOperator{\Tr}{Tr}
\def\Yc{{\cal Y}}
\def\Zc{{\cal Z}}

\begin{document}

\maketitle

\begin{abstract}
We propose several algorithms to solve McKean-Vlasov Forward Backward Stochastic Differential Equations (FBSDEs).
Our schemes rely on the approximating power of neural networks to estimate the solution or its gradient through minimization problems. As a consequence, we obtain methods able to tackle both mean-field games and mean-field control problems in moderate dimension. We analyze the numerical behavior of our algorithms on several examples including non linear quadratic models.
\end{abstract}

\vspace{5mm}


\vspace{5mm}


\section{Introduction}
\paragraph{}
This paper is dedicated to the numerical resolution in moderate dimension of the following McKean-Vlasov Forward Backward Stochastic Differential Equations (MKV FBSDEs)
\begin{equation}\label{eq: MKV FBSDE}
\begin{cases}
X_t & = \xi + \int_0^t b(s,X_s,Y_s,Z_s,\mathcal{L}(X_s),\mathcal{L}(Y_s),\mathcal{L}(Z_s))\ \di s + \int_0^t \sigma(s,X_s,\mathcal{L}(X_s))\ \di W_s \\
Y_t & = g(X_T,\mathcal{L}(X_T))  + \int_t^T f(s,X_s,Y_s,Z_s,\mathcal{L}(X_s),\mathcal{L}(Y_s),\mathcal{L}(Z_s))\ \di s - \int_t^T Z_s\ \di W_s\\
\end{cases}
\end{equation}
with $b: \R\times\R^d\times\R^k\times\R^{k\times d}\times\mathcal{P}_2(\R^d)\times\mathcal{P}_2(\R^k)\times\mathcal{P}_2(\R^{k\times d}) \mapsto \R^d$, $\sigma: \R\times\R^d \times \mathcal{P}_2(\R^d) \mapsto \R^d$, $g: \R^d\times\mathcal{P}_2(\R^d) \mapsto \R^k$, and $f: \R\times\R^d\times\R^k\times\R^{k\times d}\times\mathcal{P}_2(\R^d)\times\mathcal{P}_2(\R^k)\times\mathcal{P}_2(\R^{k\times d}) \mapsto \R^k$.  \\
$W_t$ is a $d$-dimensional $\mathcal{F}_t$-Brownian motion where $(\Omega,\mathcal{A},\mathcal{F}_t,\Prob)$ is a given filtered probability space and $T>0$.\\
$\xi$ is a given random variable in $L^2(\Omega,\mathcal{F},\Prob;\R^d)$ and
$\mathcal{P}_2(\R^n)$ stands for the space of square integrable probability measures over $\R^n$ endowed with the $2-$Wasserstein distance 
\begin{equation*}
\mathcal{W}_2(\mu,\nu) = \inf\left\{\sqrt{\E[(X-X')^2]}\ |\ X,X'\in L^2(\Omega,\mathcal{F},\Prob;\R^n),\ \mathcal{L}(X) = \mu,\ \mathcal{L}(X') = \nu  \right\}.
\end{equation*}
At last $\mathcal{L}(.)$ is a generic notation for the law of a random variable.
\paragraph{}
This kind of equation is linked to non local PDEs kwown as master equations. We refer to  \cite{carmona2018probabilistic}, chapter 4 and 5 of volume 2, for an introduction on the subject.
In \cite{chassagneux2014probabilistic}, it is shown for example that under regularity conditions, when the drift $b$ is independent of  time and the law $\mathcal{L}(Z_t)$, the driver $f$ does not depend on time and the law $\mathcal{L}(Z_t)$, $\sigma$ does not depend on time, then the resolution of equation \eqref{eq: MKV FBSDE} provides a way to estimate the solution of the equation
\begin{align}\label{eq: master}
\partial_t \mathcal{U}(t,x, \mu) + & b(x, \mathcal{U}(t,x, \mu),\partial_x \mathcal{U}(t,x, \mu) \sigma(x, \mu),\mu, \eta).\partial_x \mathcal{U}(t,x, \mu)  \nonumber \\
& + \frac{1}{2} \Tr[ \partial^2_{xx} \mathcal{U}(t,x, \mu) \sigma \sigma^\top(x, \mu)] + f(x, \mathcal{U}(t,x, \mu), \partial_x \mathcal{U}(t,x, \mu) \sigma(x, \mu),\mu, \eta)  \nonumber\\
& + \int_{\R^d} \partial_\mu \mathcal{U}(t,x, \mu)(y). b(y,\mathcal{U}(t,y, \mu),\partial_x \mathcal{U}(t,x, \mu) \sigma(x, \mu),\mu, \eta) \ \di\mu(y)  \nonumber\\
&
 + \int_{\R^d} \frac{1}{2} \Tr[ \partial_x \partial_\mu \mathcal{U}(t,x, \mu)(y) 
\sigma \sigma^\top(y, \mu)]\ \di\mu(y) =0\ \mathrm{on}\ [0,T]\times\R^d\times\Pc_2(\R^d)
\end{align} with initial condition $\mathcal{U}(0,x, \mu) = g(x,\mu)$ on $\R^d\times\Pc_2(\R^d)$, and
where  $\eta$ is a notation for the image of the probability measure $\mu$ by the mapping $x\in\R^d \longrightarrow  \mathcal{U}(t,x,\mu) \in \R^{k}$.  
Under suitable assumptions, \cite{chassagneux2014probabilistic} proves that the solution to \eqref{eq: MKV FBSDE} admits the so-called decoupling field representation $Y_s = U(s,X_s,\Lc(X_s)) $ where $(t,x,\mu)\in[0,T]\times\R^d\times\Pc_2(\R^d)\mapsto U(T-t,x,\mu) $ is a classical solution to \eqref{eq: master}. See Theorem 2.9 and equation (2.12) in \cite{chassagneux2014probabilistic}. Their results are stated in the more general of dynamics depending in the joint law of $(X_t,Y_t)$ but we state them here with the particular case of dependence in the marginal laws $\Lc(X_t),\Lc(Y_t)$ to be coherent with \eqref{eq: MKV FBSDE}. This MKV FBSDE representation is used in \cite{CCD17} to build a numerical scheme for the approximation of \eqref{eq: master} when the drift $b$ is independent of $Z$.
The equation above is non local due to the integral terms
 and the term $\partial_\mu \mathcal{U}(t,x, \mu)(y)$ stands
 for the Wasserstein derivative of $\mathcal U$ in the direction of the measure at point $(t,x, \mu)$ and evaluated at the continuous coordinate $y$. 

\paragraph{}
Equations \eqref{eq: MKV FBSDE}   appear as well  as probabilistic formulations of mean-field games or mean-field controls characterizing the value function $V$ of the game.
Mean-field games are introduced by \cite{LL06} and \cite{LL06P2} to model games with interactions between many similar players. In this theory, each player's 
dynamics and cost take into account the empirical distribution of all agents. At the limit of an infinite number of players, the search for a Nash equilibrium with close loop controls boils down to a control problem concerning a representative player whose law enters in the cost and dynamics.\\
Two probabilistic approaches based on Forward Backward Stochastic Differential Equations can be used to solve these problems:
\begin{itemize}
    \item A first approach  called the {\bf Pontryagin approach} consists as shown in \cite{CD13} in applying  the strong Pontryagin principle to these control problems. Under regularity and convexity conditions,  $Y_t$ appears to be  a stochastic representation of the gradient of the value function $V$. In this case the coefficients $b,f$ of the related MKV FBSDE \eqref{eq: MKV FBSDE} do not depend on $Z_t, \mathcal{L}(Y_t),\mathcal{L}(Z_t) $ and $k = d.$ 
    \item Another approach called the {\bf Weak approach} permits to solve the optimization problem by estimating directly $Y_t$ as the value function $V$ of the problem as shown in  \cite{CL15}. In this case the coefficients $b,f$ of the related MKV FBSDE \eqref{eq: MKV FBSDE} do not depend on $Y_t, \mathcal{L}(Y_t),\mathcal{L}(Z_t) $ and $k = 1.$ 
\end{itemize}
\paragraph{}
The numerical resolution of equations \eqref{eq: MKV FBSDE} is rather difficult since:
\begin{itemize}
    \item The dynamics are coupled through both the drift and the driver of the BSDE.
    \item The McKean-Vlasov structure of the problem requires to solve a fixed point in probability spaces. 
\end{itemize}
In the linear-quadratic setting (with quadratic cost to minimize but linear dynamics) the weak approach applied to mean-field problems gives a problem in low dimension but with quadratic coupling in $Z_t$ appearing in the backward dynamic. In contrast, the Pontryagin approach exhibits a problem in potentially high dimension (the $Z$ component is a $d\times d$ matrix in this case) but with a linear coupling in $Y_t$ which is easier to solve numerically.\\
In the case of mean-field games, only the law of $X_t$ is present in the dynamic of \eqref{eq: MKV FBSDE}.
In mean-field games of controls (also called extended mean-field games, see \cite{AK20} and the references therein), individuals interact through their controls instead of their states as in the model of trade crowding in \cite{cardaliaguet2018mean}. The law of the control thus appears in the dynamic of \eqref{eq: MKV FBSDE} and may give rise to some FBSDE depending on  the law of $Z_t$ in the weak approach or the law of $Y_t$ in the Pontryagin approach.
\paragraph{}
Existence and uniqueness of a solution to the fully coupled system \eqref{eq: MKV FBSDE} are studied by \cite{carmona2018probabilistic} when the drift does not depend on the law of $Z$. Their Theorem 4.29 gives existence of a solution under a non-degeneracy condition. However, uniqueness is a priori only expected to hold in small time, as stated in Theorem 4.24 from \cite{carmona2018probabilistic}.
In the latter we will assume that existence and uniqueness hold for the MKV FBSDE we aim to solve.

\paragraph{}
In \cite{CCD17} and \cite{AVLCDC19}, tree and grid algorithms are proposed and tested in dimension 1. It is worth mentioning that these techniques suffer from the so-called curse of dimensionality and cannot be applied when the dimension describing a player state is high (typically greater than 3 or 4). This is due to the discretization of the state space.
\paragraph{}
However, new approaches using machine learning are developed since 2017 to solve non linear parabolic PDEs through a BSDE representation.
Two kinds of methods have emerged:
\begin{itemize}
    \item The first to appear are {\bf global methods} first proposed in \cite{HJE17} to solve semi-linear PDEs. They rely on a single high-dimensional optimization problem whose resolution is difficult. It consists in the training of as many neural networks as time steps by solving in a forward way the backward representation of the PDE solution.  The $Z_t$ process is represented by a different neural network $Z^\theta_i$ with parameters $\theta$ at each date $t_i$. Instead of solving the BSDE starting from the terminal condition, the method writes it down as a forward equation and an optimization problem aiming to reach the terminal condition $g(X_T)$ by minimizing a mean squared error $\E|Y_T - g(X_T)|^2$. The approach is extended to fully nonlinear equations (nonlinear in the solution, its gradient and hessian) in \cite{beck2017machine} and the authors show that the methodology can solve some equations in high dimension.
    \cite{CWNMW19} showed that it is more effective to used a single network for all dates and besides proposed an original fixed point algorithm to solve semilinear PDEs.
    \item A second kind of {\bf local methods} first proposed in \cite{HPW19} is based on local optimization problems solved at each time step in a backward way.  Contrarily to the global method, the successive optimization problems are here in moderate dimension. Each optimization step at date $t_i$ consists in the training of only two local neural networks $Y^\theta_i, Z^\theta_i$ with parameters $\theta$. For instance, instead of solving 1 optimization problem with $N$ neural networks in the global method, \cite{HPW19} solves $N$ learning problems with 2 neural networks. Moreover the resolution is simplified by the initialization of the neural networks at time $t_i$ to their previously computed values at time $t_{i+1}$, which provides a good approximation for the current value. This strategy is inspired by the standard backward resolution of BSDE with conditional expectations from \cite{BT04} and \cite{GLW05}. 
    The methodology is extended to the much more challenging case of fully nonlinear PDEs in \cite{phawarger20} by combining it with some ideas proposed in \cite{beck2019deep}. Extensive tests performed in \cite{HPW19} show that the local method gives better results than the global one, such as \cite{phawarger20} in the case of fully nonlinear dynamics. Especially, these papers show that local methods can be used with a larger time horizon $T$ than the global method.
\end{itemize}
\paragraph{}
Machine learning techniques to solve coupled FBSDEs are investigated by several authors in \cite{HL18} and \cite{JPPZ19}, and a first method for McKean-Vlasov FBSDEs with delay is studied by \cite{FZ19} for a linear quadratic equation in dimension one. Similar and more general ideas are presented and tested in dimension one in \cite{CL19b} alongside convergence results, together with an additional method directly solving mean-field control problems by minimizing the cost without writing down optimality conditions.
The resulting algorithms proposed all rely on the global approach first initiated in \cite{HJE17}.
\paragraph{}
Our paper aims to extend these methods and to propose new ones for the resolution of McKean-Vlasov FBSDEs in moderate dimension, and go beyond one dimensional examples for which standard methods are already available (see \cite{AC10,L21}). In fact, one major advantage for the use of neural networks for solving control problems is their ability to efficiently represent high-dimensional functions without using space grids. We also study the influence of the maturity $T$ on the algorithms.
We first propose to modify the previously proposed algorithm to
stabilize its convergence. 
Our modification allows us to reduce the variance of the estimators used in the dynamic of $X_t$ and $Y_t$.
Then we propose a second algorithm relaxing the fixed point iteration algorithm by adding a neural network learning the distribution of the solution thanks to a penalization in the loss function.
At last we propose a resolution scheme based on some local resolution as in \cite{HPW19}.\paragraph{}
To simplify the presentation, we consider first order interaction models, that is to say that the dependency of the drift and cost function with respect to the laws $\mathcal{L}(X_t),\mathcal{L}(Y_t),\mathcal{L}(Z_t)$ only concerns expectations in the form
\begin{equation}
    u_t = (u_t^X,u_t^Y,u_t^Z) := (\E[\varphi_1(X_t)],  \E[\varphi_2(Y_t)], \E[\varphi_3(Z_t)]),\label{eq: def m}
\end{equation} for some continuous functions $\varphi_1,\varphi_2,\varphi_3$ with adequate domains and codomains. In this framework we can rewrite by abuse of notation 
\begin{align*}
    b(s,X_s,Y_s,Z_s,\mathcal{L}(X_s),\mathcal{L}(Y_s),\mathcal{L}(Z_s)) & = b(s,X_s,Y_s,Z_s,u_s)\\
     \sigma(s,X_s,\mathcal{L}(X_s)) & = \sigma(s,X_s,u_s^X) \\
g(X_T,\mathcal{L}(X_T)) & = g(X_T,u_T^X)\\
f(s,X_s,Y_s,Z_s,\mathcal{L}(X_s),\mathcal{L}(Y_s),\mathcal{L}(Z_s)) & = f(s,X_s,Y_s,Z_s,u_s).
\end{align*}
For instance when $\varphi_1,\varphi_2,\varphi_3$ are power functions with positive integers as exponents we recover probability distribution moments.
See Remark \ref{rem: more general law} for more general cases beyond first order interaction.
\paragraph{}
We provide multidimensional tests to show how these machine learning approaches can overcome the curse of dimensionality on some test cases first coming from a mean-field game of controls: we solve the FBSDE derived from the weak approach and the Pontryagin approach. We also consider an example arising from a non linear quadratic mean-field game.
Then we compare all the methods on some general test cases of FBSDE involving linear or quadratic dependence on the processes $X_t$, $Y_t$, $Z_t$ and on their distributions.
\paragraph{}
The structure of the paper is the following: in sections \ref{schemes} and \ref{local} we describe the proposed schemes, and in section \ref{results} we provide a numerical study of our methods in dimension 10 (except for the one-dimensional Example of Section \ref{pop}). We show that our algorithms can solve non linear-quadratic models with small maturities.
\section{Machine learning global solvers}\label{schemes}
In this section we propose three {\bf global   algorithms} based on the approach in \cite{HJE17}. 

\subsection{Algorithm principle}

\paragraph{}
We propose a generalized and refined version of the Algorithm 2 from \cite{CL19b}. We recall that a similar technique with additional networks is used in \cite{FZ19} for delayed McKean-Vlasov equations but is tested only on a one dimensional linear quadratic example. Our methods also take advantage of different expectation computation methods, introduced in section \ref{sec: expectation}. 
We present in section \ref{results} several tests in dimension 10 where the laws of $X, Y, Z$ are involved. 
\paragraph{}
We consider the Euler-Maruyama discretized FBSDE system \eqref{eq: MKV FBSDE} on a regular time grid $t_k = \frac{kT}{N}$ for $k\in\llbracket0,N\rrbracket$:
\begin{equation}\label{eq: approximated FBSDE}
\begin{cases}
X_{t_{i+1}} & = X_{t_{i}} + b\left(t_i,X_{t_{i}},Y_{t_{i}},Z_{t_{i}},u_{t_{i}}\right)\ \Delta t+ \sigma\left(t_{i},X_{t_{i}},u_{t_{i}}^X \right) \Delta W_i\\
Y_{t_{i+1}} & = Y_{t_{i}} - f\left({t_{i}},X_{t_{i}},Y_{t_{i}},Z_{t_{i}},u_{t_{i}}\right)\ \Delta t + Z_{t_{i}} \Delta W_i,
\end{cases}
\end{equation}
with terminal condition $Y_{t_N} = g(X_{t_N},u_{t_N}^X)$ and initial condition $X_0 = \xi$. We recall that $u_t$ is defined in \eqref{eq: def m}. We note $\Delta t := t_{i+1} - t_i = \frac{T}{N}$ and $(\Delta W_i)_{i=0,\cdots,N-1} := (W_{t_{i+1}} - W_{t_{i}})_{i=0,\cdots,N-1}$ the Brownian increments.
In the FBSDE theory, one requires the processes $(X_{t_{i}},Y_{t_{i}},Z_{t_{i}})$ to be $\mathcal{F}_{t_{i}}$-adapted. Therefore the backward part of the system can also be written in the conditional expectation form
\begin{equation}
    \begin{cases}
        Y_{t_i} = \E[Y_{t_{i+1}} + f\left({t_{i}},X_{t_{i}},Y_{t_{i}},Z_{t_{i}},u_{t_{i}}\right)\ \Delta t | \mathcal{F}_{t_i}]\\
       Z_{t_{i}} =  \E[Y_{t_{i+1}}\frac{\Delta W_i}{\Delta t} | \mathcal{F}_{t_i}],
    \end{cases}
\end{equation} where we see how the process $Z$ is defined. This process, specific to the stochastic case, allows the $Y$ component to be $\mathcal{F}_{t}$-adapted, even though we fix its terminal condition. It is a major difference between backward ordinary differential equations and backward stochastic differential equations.\\  We see that the whole system is coupled therefore we need to design a method allowing to solve simultaneously both equations of \eqref{eq: approximated FBSDE}.
\paragraph{}
We solve the system by the Merged Deep BSDE method introduced in \cite{CWNMW19}. $Z_{t_{i}}$ is approximated by a single feedforward neural network $\Zc_{\theta^z}(t_{i}, X_{t_{i}})$ and  $Y_0$ by a neural network $\Yc_{\theta^y}(X_{0})$ with parameters $\theta = ({\theta^y},{\theta^z})$. 
With this point of view, the discretized Brownian motion $W_t$ acts as training data in the language of machine learning, so that we can generate a training set as large as desired. 
Extensive tests conducted in \cite{CWNMW19} show that the use of a Merged network improves the training in comparison with the Deep BSDE method of \cite{HJE17}. Indeed it lowers the number of parameters to learn hence reduces the complexity of the problem. It empirically improves the accuracy of the method but also makes the training faster. That is why we focus on this architecture. It is also used by \cite{CL19b} which considers controls in the form $\big(\zeta(t_i,X_i)\big)_{i=1,\dots,N}$ for a neural network $\zeta$. The use of  recurrent networks such as Long Short Term Memory networks as in \cite{FZ19} is possible  but tests achieved in \cite{CWNMW19} seem to show that is does not bring more accuracy on Markovian problems. Other alternatives may include the GroupSort network \cite{ALG19} for a better control of the  Lipschitz constant of the approximation, or some special networks preserving some properties of the solution but they have not been tested.

 The motivation for such an approximation comes from the notion of decoupling field, also used for numerical purposes in \cite{AVLCDC19} or \cite{CCD17}, which gives the existence of functions $u,v$ (see the paragraph below \eqref{eq: master}) such that
\begin{equation}
Y_t = u(t,X_t,\mathcal{L}(X_t)),\ Z_t = v(t,X_t,\mathcal{L}(X_t)).
\end{equation}
Numerically, it is enough to consider $Y$ and $Z$ as a function of the couple $(t,X_t)$. In fact, the law of the solution (and therefore its moments) can be seen as a function of $t$. That's why we search for a representation 
\begin{equation}
Y_t = \widetilde{u}(t,X_t),\ Z_t = \widetilde{v}(t,X_t).
\end{equation} The forward-backward system is transformed into a forward system and an optimization problem aiming to satisfy the terminal condition of the BSDE through the loss function $\E\big[\big(Y_T - g\big(X_T,$ $u^X_{t_N}\big)
 \big)^2\big]$. To simplify notations, $X_i := X_{t_i}$ and similarly for $Y$ and $Z$.

\paragraph{}
In practice the loss function is minimized with the Adam gradient descent method \cite{KB14}. 
In any case, the goal of our scheme is to learn both the optimal control and the distribution of $X_t, Y_t, Z_t$. In the following, $B$ is the batch size, that is the number of particles we will sample to approximate the loss function in expectation form, see for instance \eqref{eq: batch size}. $N$ is the number of time steps and $M$ is the number of previous batches expectations to keep in memory. 
\paragraph{}
We use for $\Zc$ and $\Yc$ feedforward neural networks with 3 hidden layers ($d + 10$ neurons in each) with hyperbolic tangent function as activation functions and an output layer with identity as activation. It is worth noticing that because the merged neural network takes the couple $(t,X)$ as inputs, we cannot use batch normalization since the distribution of $X_i$ is not stationary over time.

\subsection{Estimation of the expectation}\label{sec: expectation}
\paragraph{}
A key step for the methods is to estimate the mean-field parameter $u$. It has a significant effect on the algorithms performances. We note $\theta_m=({\theta^y_m},{\theta^z_m})$ the neural network parameters at optimization iteration $m$ and $u_i = (u_i^X,u_i^Y,u_i^Z)$ the estimation of 
$u_{t_i}$. 
In the algorithms described below, the approximated processes are considered as functions of the parameters $\theta$ of the neural network.\\
Several methods can be used to approximate the moments of the solution involved in the stochastic McKean-Vlasov dynamics:
\begin{itemize}
    \item \textbf{Direct}: use the empirical mean of the current batch of particles
    \begin{equation} \label{eq: Monte-Carlo}
       u_i = \frac{1}{B} \left(\sum_{j=1}^{B} \varphi_1(X^j_{i}(\theta_m)), \sum_{j=1}^{B} \varphi_2(Y^j_{i}(\theta_m)), \sum_{j=1}^{B} \varphi_3(Z^j_{i}(\theta_m))\right),\ i=0, \cdots , N-1.
    \end{equation}
    
    Alternatively one could use instead the last batch particles
    to estimate the law 
    \begin{equation} 
       u_i = \frac{1}{B} \left(\sum_{j=1}^{B} \varphi_1(X^j_{i}(\theta_{m-1})), \sum_{j=1}^{B} \varphi_2(Y^j_{i}(\theta_{m-1})), \sum_{j=1}^{B} \varphi_3(Z^j_{i}(\theta_{m-1}))\right),\ i=0, \cdots , N-1.
    \end{equation} The difference lies in the fact that in one case the optimization of the parameters $\theta_m$ at iteration $m$ modifies the current estimation of the law whereas using the previously computed parameters $\theta_{m-1}$ fixes the law and simplifies the optimization problem.
     In practice, for the numerical tests of Section \ref{results} we use the formula \eqref{eq: Monte-Carlo}. This approach requires to handle very large batches, typically of the order of $B=10,000$ sample paths get a reasonable approximation of the laws. This is the approach used by \cite{FZ19} and \cite{CL19b}. 
    
    We solve the following optimization problem
    \begin{align}
        \min_{\theta = (\theta^y,\theta^z)}\ & \frac{1}{B} \sum_{k=1}^B \Big|Y_N^k(\theta) - g\Big(X_N^k(\theta),\frac{1}{B} \sum_{j=1}^B \varphi_1(X_N^j(\theta))\Big) \Big|^2\label{eq: batch size}\\
         X_{i+1}^j(\theta) =\ & X_i^j(\theta) + b\left(t_i,X_i^j(\theta),Y_i^j(\theta),\Zc_{\theta^z}\left(t_i,X_i^j(\theta)\right) ,u_{i} \right) \Delta t 
         +  \sigma\left(t_i,X_i^j(\theta) ,u_{i}^{X} \right) \Delta W_i^j\nonumber\\
         Y_{i+1}^j(\theta) =\ & Y_i^j(\theta) - f\left({t_{i}},X_{i}^j(\theta),Y_i^j(\theta),\Zc_{\theta^z}\left(t_i,X_{i}^j (\theta)\right) ,u_{i} \right)\ \Delta t 
         +  \Zc_{\theta^z}\left(t_i,X_{i}^j(\theta)\right) \Delta W_i^j \nonumber\\
         u_i =\ & \frac{1}{B} \left(\sum_{j=1}^{B} \varphi_1(X^j_{i}(\theta)), \sum_{j=1}^{B} \varphi_2(Y^j_{i}(\theta)), \sum_{j=1}^{B} \varphi_3\left(\Zc_{\theta^z}\left(t_i,X_{i}^j(\theta)\right)\right)\right)\nonumber\\
         X_0^j =\ & \xi^j \sim \xi,\ j=1,\cdots,B \nonumber\\
         Y_0^j(\theta) =\ & \Yc_{\theta^y}(X_0^j),\ \nonumber  %
         \\          i = \ & 1,\cdots,N-1.\nonumber
    \end{align}
    
    The {\bf Global Direct solver} leads to algorithm \ref{algo: direct}. 
    
\begin{algorithm}[H]
\caption{Global Direct solver }\label{algo: direct}
\begin{algorithmic}[1]
\State Let $\Yc_{\theta^y}(\cdot)$ be a neural network with parameter  ${\theta^y}$, defined on $\R^d$ and valued in $\R^{k}$ 
, $\Zc_{\theta^z}(\cdot,\cdot)$ be a neural network with parameter  ${\theta^z}$, defined on $\R^{+} \times \R^d$ and valued in $\R^{k\times d}$, so that $ \theta = ({\theta^y}, {\theta^z})$ is initialized with value $\theta_0 = (\theta^y_0, \theta^z_0)$.
\For{$m$ from 0 to $K$} \Comment{Stochastic gradient iterations}
\State Sample $(\xi^j)_{j=1,\cdots,B}$ from $B$ independent copies of the initial condition $\xi$.
\State Set $\forall j\in\llbracket1,B\rrbracket,\ X_0^j(\theta_m) = \xi^j \in\R^d , Y_0^j(\theta_m) = \Yc_{\theta^y_m}(\xi^j) \in\R^k.$
\For{$i$ from 0 to $N-1$}
\State $u_{i} = (u_{i}^X,u_{i}^Y,u_{i}^Z) = \frac{1}{B}\sum_{j=1}^{B} \left( \varphi_1(X^j_{i}(\theta_m)),  \varphi_2(Y^j_{i}(\theta_m)),  \varphi_3\Big(\Zc_{\theta^z_m}\left(t_i,X_{i}^j(\theta_{m})\right)\Big)\right)$
\For{$j$ from 1 to $B$}
\State Sample $\delta_i^j$ from a $d$-dimensional standard Gaussian vector.
\State $X_{i+1}^j(\theta_m) = X_i^j(\theta_m) + b\left(t_i,X_i^j(\theta_m),Y_i^j(\theta_m),\Zc_{\theta^z_m}\left(t_i,X_i^j(\theta_m)\right) ,u_{i}\right) \Delta t + \sqrt{\Delta t}\ \sigma\left(t_i,X_i^j(\theta_m) ,u_{i}^X\right) \delta_i^j$
\State $ Y_{i+1}^j(\theta_m) = Y_i^j(\theta_m) - f\left({t_{i}},X_{i}^j(\theta_m),Y_i^j(\theta_m),\Zc_{\theta^z_m}\left(t_i,X_{i}^j (\theta_m)\right) ,u_{i}\right)\ \Delta t +\sqrt{\Delta t}\ \Zc_{\theta^z_m}\left(t_i,X_{i}^j(\theta_m)\right) \delta_i^j$
\EndFor
\EndFor
\State $\overline{X_{N}}(\theta_m) = \frac{1}{B} \sum_{j=1}^{B} \varphi_1(X^j_{N}(\theta_m))$,
\State $J(\theta_m) = \frac{1}{B}\sum_{j=1}^{B} \left(Y_N^j(\theta_m) - g\left(X_N^j(\theta_m),\overline{X_{N}}(\theta_m)\right) \right)^2$
\State Calculate $\nabla J(\theta_m)$ by back-propagation.
\State Update $\theta_{m+1} = \theta_m - \rho_m \nabla J(\theta_m)$.
\EndFor
\end{algorithmic}
\end{algorithm}
    \item \textbf{Dynamic}: a method which dynamically updates the estimation on $(M+1) B$ samples. The expectations from the last $M$ batches are kept in memory in an array $$(\zeta_{i,r})_{\begin{array}{c}
    i=0, \dots, N-1, \\
    r= 0, \dots ,  M-1
    \end{array} } 
    $$ 
    initialized with values $(\E[\varphi_1(\xi)], \varphi_2(0),\varphi_3(0))^{N\times M}$.\\
    At iteration $m - 1$, $\nu_{{i}}^{(m-1)}$ is defined as the empirical mean on these previous sample paths. On a new batch, the expectation is computed by averaging the previous estimation $\nu_{{i}}^{(m-1)}$ and the current batch empirical mean  by the following algorithm  used for $i=0, \cdots , N-1$:
    \begin{align} \label{eq: dynamic Monte-Carlo}
            \nu_{i}^{(m-1)} &= \frac{1}{M} \sum_{r=0}^{M-1} \zeta_{i,r}, \nonumber\\
         u_i &= \frac{M \nu_{{i}}^{(m-1)} + \frac{1}{B}\left(\sum_{j=1}^{B} \varphi_1(X^j_{i}(\theta_m)),\sum_{j=1}^{B} \varphi_2(Y^j_{i}(\theta_m)), \sum_{j=1}^{B} \varphi_3(Z^j_{i}(\theta_m))\right)}{M + 1}, \nonumber \\
        \zeta_{i,m\% M} & = \frac{1}{B}\left(\sum_{j=1}^{B} \varphi_1(X^j_{i}(\theta_m)),\sum_{j=1}^{B} \varphi_2(Y^j_{i}(\theta_m)), \sum_{j=1}^{B} \varphi_3(Z^j_{i}(\theta_m))\right).
    \end{align} The notation $m\%M$ refers to the remainder of the Euclidian division of $m$ by $M$. This technique allows to use smaller batches of size 100 or 1000. Thus it is more efficient in terms of convergence speed in comparison with the direct approach. This method can be seen as a dynamic fixed point approach. 
    \paragraph{}
    The idea behind this update rule comes from online learning in machine learning. $1/M$ can be interpreted as a learning rate quantifying the updating speed. From the current estimation of the particles law, we introduce a small correction related to the new observed samples. Therefore the estimation is much more stable through iterations compared to the instantaneous update of the law used by the Direct method. After $M$ batches, the older samples are forgotten, since they don't represent anymore the current law. Indeed we expect the convergence for a good choice of $M$. If this parameter is too small the stabilization would be inefficient and on the contrary a too large $M$ would slow down the learning process by introducing a bias in the law.
    For instance in our numerical experiments of Section \ref{results} we use $M = 100$ for a total of 2000 gradient descent iterations.
    
    \begin{Rem}\label{rem: more general law}
    	If the law dependence is more general than a first order interaction 
    	and is given by a continuous function $F: \mu \in\mathcal{P}_2(\R^d) \mapsto \R^k$ then the Direct method can be straightforwardly applied to the equation by estimating $F(\mathcal{L}(X_t))$ by the so-called empirical projection $ F(\frac{1}{B} \sum_{j=1}^B \delta_{X^j_{t_i}})$ for identically distributed particles $(X^j_{t_i})_{j=1,\dots,B}$ on a time grid $t_0,\cdots,t_N$. Concerning the Dynamic approach, it would require to keep in memory the previously computed particles from the last M batches which is costly.
    \end{Rem}
    
    \begin{Rem}
    The fixed point approach is known to be convergent theoretically only for small maturities. In practice, the theoretical bound on the maturity found on the simple example given for example in paragraph 3.1 in \cite{AVLCDC19} is far too pessimistic.  We will see that the restriction is not relevant on all our test cases.
    \end{Rem}

    For a given iteration $m$, given the estimations $(\zeta_{i,r})_{\begin{array}{c}
    i=0, \dots, N-1, \\
    r= 0, \dots ,  M-1
    \end{array} }$ of $u_i$ on the last $M$ iterations, we perform one gradient descent step for the following optimization problem
    \begin{align*}
        \min_{\theta_m=({\theta^y_m},{\theta^z_m})}\ & \frac{1}{B} \sum_{k=1}^B \Big|Y_N^k(\theta_m) - g\Big(X_N^k(\theta_m),\frac{1}{B} \sum_{j=1}^B \varphi_1(X_N^j(\theta_m))\Big) \Big|^2\\
         X_{i+1}^j(\theta_m) =\ & X_i^j(\theta_m) + b\left(t_i,X_i^j(\theta_m),Y_i^j(\theta_m),\Zc_{{\theta^z_m}}\left(t_i,X_i^j(\theta_m)\right),\widetilde{u_{i}}\right) \Delta t \\ & +  \sigma\left(t_i,X_i^j(\theta_m),\widetilde{u_{i}}^X(\theta_m)\right) \Delta W_i^j\\
         Y_{i+1}^j(\theta_m) =\ & Y_i^j(\theta_m) - f\left({t_{i}},X_{i}^j(\theta_m),Y_i^j(\theta_m),\Zc_{{\theta^z_m}}\left(t_i,X_{i}^j (\theta_m)\right),\widetilde{u_{i}}\right)\ \Delta t \\ & +  \Zc_{{\theta^z_m}}\left(t_i,X_{i}^j(\theta_m)\right) \Delta W_i^j \\
         u_i =\ & \frac{1}{B} \left(\sum_{j=1}^{B} \varphi_1(X^j_{i}(\theta_m)), \sum_{j=1}^{B} \varphi_2(Y^j_{i}(\theta_m)), \sum_{j=1}^{B} \varphi_3\left(\Zc_{{\theta^z_m}}\left(t_i,X_{i}^j(\theta_m)\right)\right)\right)\\
         \widetilde{u_{i}} =\ & \frac{\sum_{r=0}^{M-1} \zeta_{i,r} + u_{i} }{M + 1}\\
         X_0^j =\ & \xi^j \sim \xi,\ j=1,\cdots,B\\
         Y_0^j(\theta_m) =\ & \Yc_{\theta^y_m}(X_0^j)\\
         i = \ & 1,\cdots,N-1.
    \end{align*} Then we update $(\xi_{i})_i$ by forgetting the oldest estimation and keeping in memory the new one,  $(u_i)_i$ (see \eqref{eq: dynamic Monte-Carlo}). The {\bf Global Dynamic solver} is given more explicitly in algorithm \ref{algo: dynamic}.
\begin{algorithm}[H]
\caption{Global Dynamic solver}\label{algo: dynamic}
\begin{algorithmic}[1]
\State Let $\Yc_{\theta^y}(\cdot)$ be a neural network with parameter  ${\theta^y}$, defined on $\R^d$ and valued in $\R^{k}$
, $\Zc_{\theta^z}(\cdot,\cdot)$ be a neural network with parameter  ${\theta^z}$, defined on $\R^{+} \times \R^d$ and valued  in $\R^{k\times d}$, so that $ \theta = ({\theta^y}, {\theta^z})$ is initialized with value $\theta_0=(\theta^y_0, \theta^z_0)$.
\State Set $\forall i\in\llbracket0,N-1\rrbracket,\ \forall r\in\llbracket0,M-1\rrbracket,\ \zeta_{i,r} = (\E[\varphi_1(\xi)], \varphi_2(0), \varphi_3(0))$.
\For{$m$ from 0 to $K$}
\State Sample $(\xi^j)_{j=1,\cdots,B}$ from $B$ independent copies of the initial condition $\xi$.
\State Set $\forall j\in\llbracket1,B\rrbracket,\ X_0^j(\theta_m) = \xi^j\in\R^d, Y_0^j(\theta_m) = \Yc_{\theta^y_m}(\xi^j)\in\R^k.$
\For{$i$ from 0 to $N-1$}
\State $u_{i} = \frac{1}{B} \left(\sum_{j=1}^{B} \varphi_1(X^j_{i}(\theta_m)),  \sum_{j=1}^{B} \varphi_2(Y^j_{i}(\theta_m)),  \sum_{j=1}^{B} \varphi_3\Big(\Zc_{\theta^z_m}\left(t_i,X_{i}^j(\theta_{m})\right)\Big)\right)$
\State $\widetilde{u_{i}} = (\widetilde{u_{i}}^X,\widetilde{u_{i}}^Y,\widetilde{u_{i}}^Z) =\frac{\sum_{r=0}^{M-1} \zeta_{i,r} + u_{i}}{M + 1}$
\For{$j$ from 1 to $B$}
\State Sample $\delta_i^j$ from a $d$-dimensional standard Gaussian vector.
\State $X_{i+1}^j(\theta_m) = X_i^j(\theta_m) + b\left(t_i,X_i^j(\theta_m),Y_i^j(\theta_m),\Zc_{\theta^z_m}\left(t_i,X_i^j(\theta_m)\right),\widetilde{u_{i}} \right) \Delta t + \sqrt{\Delta t}\ \sigma\left(t_i,X_i^j(\theta_m),\widetilde{u_{i}}^{X} \right) \delta_i^j$
\State $Y_{i+1}^j(\theta_m) = Y_i^j(\theta_m) - f\left({t_{i}},X_{i}^j(\theta_m),Y_i^j(\theta_m),\Zc_{\theta^z_m}\left(t_i,X_{i}^j (\theta_m)\right),\widetilde{u_{i}} \right)\ \Delta t + \sqrt{\Delta t}\ \Zc_{\theta^z_m}\left(t_i,X_{i}^j(\theta_m)\right) \delta_i^j$
\EndFor
\State $\zeta_{i,m\% M} = u_{i} $
\EndFor
\State $\overline{X_{N}}(\theta_m) = \frac{1}{B} \sum_{j=1}^{B} \varphi_1(X^j_{N}(\theta_m))$,
\State $J(\theta_m) = \frac{1}{B}\sum_{j=1}^{B} \left(Y_N^j(\theta_m) - g\left(X_N^j(\theta_m),\overline{X_{N}}(\theta_m)\right) \right)^2$
\State Calculate $\nabla J(\theta_m)$ by back-propagation.
\State Update $\theta_{m+1} = \theta_m - \rho_m \nabla J(\theta_m)$.
\EndFor
\end{algorithmic}
\end{algorithm}

    \item \textbf{Expectation}: estimate $u_t$ by a neural network $\Psi_{\theta^\Psi}$ with input $t$ and parameters ${\theta^\Psi}$. 
    \begin{equation} \label{eq: neural expectation}
        u_i (\theta^\Psi) = \Psi_{\theta^\Psi}(t_i) = ( \Psi_{\theta^\Psi}^X(t_i),\Psi_{\theta^\Psi}^Y(t_i) ,\Psi_{\theta^\Psi}^Z(t_i) ),\ i=0, \cdots , N.
    \end{equation}
    A penalization term
    \begin{equation*}
    \E\left[\frac{\lambda}{N}\sum_{i=0}^{N-1} \left\lVert\Psi_{\theta^\Psi}(t_i) - \frac{1}{B}\left( \sum_{j=1}^{B} \varphi_1(X^j_{i}(\theta)), \sum_{j=1}^{B} \varphi_2(Y^j_{i}(\theta)), \sum_{j=1}^{B} \varphi_3(Z^j_{i}(\theta))\right)\right\rVert_2^2\right],
    \end{equation*}
    is added to the loss function. 
    We will see that in practice this method is quite involved to use because the performances heavily depend upon the choice of the parameter $\lambda$. This approach provides a relaxation of the fixed point method.
    
    We solve the following optimization problem
    \begin{align*}
        \min_{\theta = (\theta^y,\theta^z,\theta^\Psi)}\ & \frac{1}{B} \sum_{k=1}^B \Big|Y_N^k(\theta) - g\Big(X_N^k(\theta),\frac{1}{B} \sum_{j=1}^B \varphi_1(X_N^j(\theta))\Big) \Big|^2\\
        & + \frac{\lambda}{N}\sum_{i=0}^{N-1}\left\lVert u_i(\theta^\Psi) - \Psi_{{\theta^\Psi}}(t_i)\right\rVert^2\\
         X_{i+1}^j(\theta) =\ & X_i^j(\theta) + b\left(t_i,X_i^j(\theta),Y_i^j(\theta),\Zc_{\theta^z}\left(t_i,X_i^j(\theta)\right),\Psi_{\theta^\Psi}(t_i)\right) \Delta t \\ & +  \sigma\left(t_i,X_i^j(\theta),\Psi_{\theta^\Psi}(t_i)^X\right) \Delta W_i^j\\
         Y_{i+1}^j(\theta) =\ & Y_i^j(\theta) - f\left({t_{i}},X_{i}^j(\theta),Y_i^j(\theta),\Zc_{\theta^z}\left(t_i,X_{i}^j (\theta)\right),\Psi_{\theta^\Psi}(t_i)\right)\ \Delta t \\ & +  \Zc_{\theta^z}\left(t_i,X_{i}^j(\theta)\right) \Delta W_i^j \\
         u_i =\ & \frac{1}{B} \left(\sum_{j=1}^{B} \varphi_1(X^j_{i}(\theta)), \sum_{j=1}^{B} \varphi_2(Y^j_{i}(\theta)), \sum_{j=1}^{B} \varphi_3\left(\Zc_{\theta^z}\left(t_i,X_{i}^j(\theta)\right)\right)\right)\\
         X_0^j =\ & \xi^j \sim \xi,\ j=1,\cdots,B\\
         Y_0^j(\theta) =\ & \Yc_{\theta^y}(X_0^j)\\
         i = \ & 1,\cdots,N-1.
    \end{align*}
    The {\bf Global Expectation solver} is described in algorithm \ref{algo: expectation}. 
    The parameter $\lambda$ is chosen by trial and error.
    \begin{algorithm}[H]
\caption{Global Expectation solver}\label{algo: expectation}
\begin{algorithmic}[1]
\State Let $\Yc_{\theta^y}(\cdot)$ be a neural network with parameter  ${\theta^y}$, defined on $\R^d$ and valued in $\R^{k}$, $\Zc_{\theta^z}(\cdot,\cdot)$ defined on $\R^{+} \times \R^d$, $\Psi_{\theta^\Psi}(\cdot)=( \Psi_{\theta^\Psi}^X(\cdot),\Psi_{\theta^\Psi}^Y(\cdot) ,\Psi_{\theta^\Psi}^Z(\cdot) )$ defined on $\R^{+}$ be neural networks with parameters $\theta^z$, ${\theta^\Psi}$, taking values respectively in $\R^{k\times d}$ and $\R^{d}\times\R^{k}\times\R^{k\times d}$, so that $\theta = ({\theta^y}, \theta^z, \theta^\Psi)$ is initialized with value $\theta_0 = (\theta^y_0, \theta^z_0, \theta^\Psi_0)$.
\For{$m$ from 0 to $K$}
\State Sample $(\xi^j)_{j=1,\cdots,B}$ from $B$ independent copies of the initial condition $\xi$.
\State Set $\forall j\in\llbracket1,B\rrbracket,\ X_0^j(\theta_m) = \xi^j\in\R^d, Y_0^j(\theta_m) = \Yc_{\theta^y_m}(\xi^j)\in\R^k$.
\For{$i$ from 0 to $N-1$}
\State $u_{i} = \frac{1}{B}\left( \sum_{j=1}^{B} \varphi_1(X^j_{i}(\theta_m)), \sum_{j=1}^{B} \varphi_2(Y^j_{i}(\theta_m)),\sum_{j=1}^{B} \varphi_3\Big(\Zc_{\theta^z_m}\left(t_i,X_{i}^j(\theta_m)\right)\Big)\right) $
\For{$j$ from 1 to $B$}
\State Sample $\delta_i^j$ from a $d$-dimensional Gaussian vector.
\State $X_{i+1}^j(\theta_m) = X_i^j(\theta_m) + b\left(t_i,X_i^j(\theta_m),Y_i^j(\theta_m),\Zc_{\theta^z_m}\left(t_i,X_i^j(\theta_m)\right),\Psi_{\theta^\Psi_m}(t_i)\right) \Delta t + \sqrt{\Delta t}\ \sigma\left(t_i,X_i^j(\theta_m),\Psi_{\theta^\Psi_m}^X(t_i)\right) \delta_i^j$
\State $Y_{i+1}^j(\theta_m) = Y_i^j(\theta_m) - f\left({t_{i}},X_{i}^j(\theta_m),Y_i^j(\theta_m),\Zc_{\theta^z_m}\left(t_i,X_{i}^j (\theta_m)\right),\Psi_{\theta^\Psi_m}(t_i)\right)\ \Delta t + \sqrt{\Delta t}\ \Zc_{\theta^z_m}\left(t_i,X_{i}^j(\theta_m)\right) \delta_i^j$
\EndFor
\EndFor
\State $\overline{X_{N}}(\theta_m) = \frac{1}{B} \sum_{j=1}^{B} \varphi_1(X^j_{N}(\theta_m))$,
\State $J(\theta_m) = \frac{1}{B}\sum_{j=1}^{B} \left(Y_N^j(\theta_m) - g\left(X_N^j(\theta_m), \overline{X_{N}}(\theta_m)\right) \right)^2 + \frac{\lambda}{N}\sum_{i=0}^{N-1}\left(u_{i} - \Psi_{\theta^\Psi_m}(t_i)\right)^2$
\State Calculate $\nabla J(\theta_m)$ by back-propagation.
\State Update $\theta_{m+1} =\theta_m - \rho_m \nabla J(\theta_m)$.
\EndFor
\end{algorithmic}
\end{algorithm}
\end{itemize}
We will compare the performances of these techniques on several examples in section \ref{results}.

\section{A local solver}\label{local}

We also propose a local method inspired by the Deep Backward Dynamic Programming introduced by \cite{HPW19} and \cite{phawarger20}. It considers local minimization problems between contiguous time steps. In this case there are as many networks as time steps. We replace a global optimization setting by a set of smaller problems.
\paragraph{}
In this method for $i\in\llbracket 0, N-1\rrbracket$, $Z_{{i}}$ and $Y_{{i}}$ are approximated by a neural network $(\Zc^i_{\theta_{i}^z}(\cdot)$ $,\Yc^i_{\theta_{i}^y}(\cdot))$ with parameters $\theta = (\theta_{0}^y,\theta_{0}^z,\cdots,\theta_{N-1}^y, \theta_{N-1}^z)$. At iteration $m$, with $\theta_m = (\theta_{m,0}^y , \theta_{m,0}^z,\cdots, \theta_{m,N-1}^y,$ $\theta_{m,N-1}^z)$ , we simulate $X_i(\theta_m)$ with the previously computed parameters $\theta_m$. 
\begin{align}\label{eq: dynamics X local}
    X_{i+1}^j(\theta_{m}) =\ & X_i^j(\theta_{m}) + b\bigg(t_i,X_i^j(\theta_{m}),\Yc^i_{\theta_{m,i}^y}\left(X_i^j(\theta_{m})\right),\Zc^i_{\theta_{m,i}^z}\left(X_i^j(\theta_{m})\right),\widetilde{u_{i}}\bigg)\ \Delta t \\ +\ &   \sigma\left(t_{i},X_i^j(\theta_{m}),\widetilde{u_{i}}^{X}\right) \Delta W_i^j.\nonumber
\end{align}
This first step allows to find the areas visited by the controlled process. Using $R$ samples, we compute the empirical mean $\overline{m_i^m} = \frac{1}{R} \left( \sum_{j=1}^{R} X^j_{i}(\theta_{m})\right)$ and variance $V_i^m = \frac{1}{R} \sum_{j=1}^{R} \left(X^j_{i}(\theta_{m})\right)^2 - \frac{1}{R^2} \left(\sum_{j=1}^{R} X^j_{i}(\theta_{m})\right)^2$ of $X_i(\theta_m)$. We estimate $u_t$ as in the Dynamic method \eqref{eq: dynamic Monte-Carlo}: 
\begin{align*}
    u_i =\ & \frac{1}{R} \left(\sum_{j=1}^{R} \varphi_1(X^j_{i}(\theta_{m})), \sum_{j=1}^{R} \varphi_2\left(\Yc^i_{\theta_{m,i}^y}\left(t_i,X_{i}^j(\theta_{m})\right)\right), \sum_{j=1}^{R} \varphi_3\left(\Zc^i_{\theta_{m,i}^z}\left(t_i,X_{i}^j(\theta_{m})\right)\right)\right)\\
         \widetilde{u_{i}}  =\ & \frac{\sum_{r=0}^{M-1} \zeta_{i,r} + u_{i} }{M + 1}.
\end{align*}
A priori $R$ can be different from the batch size $B$. It is the batch size used for the estimation of the law with the previously computed parameters. In the numerical tests of Section \ref{results} we use $B=100$ or $B=300$ for the backward optimization and a larger value $R= 50000$ for the Monte-Carlo forward estimation of the law. Then we solve backward problems to find the $\theta_{m+1}$ by sampling $B$ independent copies of $X_{i}$ through a Gaussian distribution $\mathcal{N}(\overline{ m_i^m},V_i^m )$ with frozen parameters $\theta_m$ : 
\begin{itemize}
    \item First sample $B$ independent copies $X^1_N,\cdots,X^B_N$ of $X_{N}$ following a Gaussian distribution $\mathcal{N}(\overline{ m_N^m},V_N^m)$. $Y^N_{\theta_{m+1,N}^y}(X_N^j)$ is set to the terminal condition $g\left(X_N^j ,u_N^X\right)$.
\item For $i$ from $N-1$ to $0$:\\
Sample $B$ independent copies $X^1_i,\cdots,X^B_i$ of $X_{i}$ following a Gaussian distribution $\mathcal{N}(\overline{m_i^m},$ $V_i^m)$. Diffuse according to the dynamics \eqref{eq: dynamics X local} of $X$, starting from $X^1_i,\cdots,X^B_i$, to obtain $X^1_{i+1},\cdots,X^B_{i+1}$.\\
Solve the local optimization problem 
\begin{align*}
\min_{\theta=(\theta^y,\theta^z)} \frac{1}{B} \sum_{j=1}^B & \left|Y^{i+1}_{\theta_{m+1,i+1}^y}\left(X_{i+1}^j\right) - {Y^i_{\theta^y}\left(X_i^j\right)}+ f\left({t_{i}},X_i^j,{\Yc^i_{\theta^y}\left(X_i^j\right)},{\Zc^i_{\theta^z}\left(X_i^j\right)},\widetilde{u_{i}}\right) \Delta t \right. \\ & \left. - {\Zc^i_{\theta^z}\left(X_i^j\right)} \Delta W_{i}^j \right|^2,
\end{align*}starting from the parameter value $\theta_{m+1,i+1}$. We can then update the $\theta$ value by denoting as  $\theta_{m+1,i}$ the argmin value of the minimization problem.
\item Repeat the previous steps for the iteration $m+1$ until reaching $K$ iterations.
\end{itemize}

In the version of the {\bf Local Dynamic solver} given in algorithm \ref{algo: local}, we use the dynamic update of the expectations introduced previously in the dynamic solver of section \ref{schemes}. In this algorithm $H$ stands for the number of gradient steps to perform at each step of the algorithm and  $R$ is the number of samples for the laws estimation.

\begin{Rem}
Because we have to learn the dynamic of the forward process, the use of a backward resolution is not as obvious as in \cite{HPBL18,BHPL19}. We have to alternate between forward dynamic estimations and backward resolutions. \\
More precisely, here we solve fully coupled FBSDEs, whereas the works \cite{HPW19,phawarger20} consider decoupled FBSDEs, where the forward process $X$ can be simulated independently of $Y,Z$. Estimating the law of $X$ and sampling from a normal distribution allows us to decouple and solve locally the FBSDEs in areas visited by $X$ and its law. However, because of this freezing of the forward dynamics, another fixed point problem has to be solved. Other approaches for fully coupled FBSDEs like \cite{HL18} and \cite{JPPZ19} rely on the global machine learning method initiated by \cite{HJE17}.\\
Notice that we estimate the law dependent part of the dynamics in the first loop but for the resolution loop, we use a normal distribution as training measure. 
Another natural training measure could be the empirical law of the process. Here, we employ an hybrid method between the two methods proposed in Remark 2.1 from \cite{BHPL19}. We exploit the approximated control to estimate the law of $X$ and use a Gaussian law with the same two first moments as $X$ as training measure. In the numerical examples of \cite{BHPL19}, the authors also use Gaussian distributions as training measures, but with fixed parameters through the iterations. Here we adapt the training measures during the training to have the same expectation and variance as the controlled process thanks to an exploration step.
\end{Rem}

\begin{algorithm}[H]
\caption{Local Dynamic solver}\label{algo: local}
\begin{algorithmic}[1]
\State Let $(\Yc^i_{\theta^y_i}(\cdot), \Zc^i_{\theta^z_i}(\cdot))$ be some neural networks defined on $\R^d$ with values in $\R^k\times\R^{k\times d}$ for $i=0, \cdots, N-1$ and parameters $\theta = (\theta^y_0,\theta^z_0,\cdots, \theta^y_{N-1},\theta^z_{N-1})$ initialized with values $\theta_0 =(\theta^y_{0,0},\theta^z_{0,0},\cdots, \theta^y_{0,N-1},\theta^z_{0,N-1}).$
\State Set $\forall i\in\llbracket0,N\rrbracket,\ \forall r\in\llbracket0,M-1\rrbracket,\ \zeta_{i,r} = (\E[\xi], 0, 0)$.
\For{$m$ from 0 to $K$}
\State Sample $\delta_i^j$ from a $d$-dimensional standard Gaussian vector,  $i =0, \cdots, N$, $j=1, \cdots, R$.
\State Sample $(\xi^j)_{j=1,\cdots,B}$ from $R$ independent copies of the initial condition $\xi$.
\State Set $\forall j\in\llbracket1,R\rrbracket,\ X_0^j(\theta_m) = \xi^j \in\R^d.$
\For{$i$ from 0 to $N$} \Comment{Forward estimation of the laws}
\State $l_{i}^m = \frac{1}{R} \left( \sum_{j=1}^{R} X^j_{i}(\theta_{m})\right)$
\State $u_{i}  = \frac{1}{R} \left( \sum_{j=1}^{R} \varphi_1(X^j_{i}(\theta_{m})), \sum_{j=1}^{R} \varphi_2(\Yc^i_{\theta_{m,i}^y}\left(X_i^j(\theta_{m})\right)),
\sum_{j=1}^{R} \varphi_3\Big(\Zc^i_{\theta_{m,i}^z}\left(X_i^j(\theta_{m})\right)\Big) \right)$
\State $V_{i}^m = \frac{1}{R} \sum_{j=1}^{R} \left(X^j_{i}(\theta_{m})\right)^2 - \frac{1}{R^2} \left(\sum_{j=1}^{R} X^j_{i}(\theta_{m})\right)^2$
\State $\widetilde{u_{i}} = (\widetilde{u_{i}}^X,\widetilde{u_{i}}^Y,\widetilde{u_{i}}^Z)  =\frac{\sum_{r=0}^{M-1} \zeta_{i,r} + u_{i} }{M + 1}$
\State $\zeta_{i,m\% M} = u_{i}( \theta_m)$
\For{$j$ from 1 to $R$}
\State $X_{i+1}^j(\theta_{m}) = X_i^j(\theta_{m}) + b\bigg(t_i,X_i^j(\theta_{m}),\Yc^i_{\theta_{m,i}^y}\left(X_i^j(\theta_{m})\right),\Zc^i_{\theta_{m,i}^z}\left(X_i^j(\theta_{m})\right),\widetilde{u_{i}}\bigg)\ \Delta t + \sqrt{\Delta t}\ \sigma\left(t_{i},X_i^j(\theta_{m}),\widetilde{u_{i}}^{X} \right) \delta_i^j$
\EndFor
\EndFor
\For{$i$ from $N-1$ to 0} \Comment{Backward resolution}
\State $\hat \theta_0= \theta_{m,i}$ 
\For{$h$ from 0 to $H-1$} \Comment{Gradient descent with simulated data for $X$}
\For{$j$ from 1 to $B$}
\State Sample $\Xi_i^j, \Theta_i^j$ from $d$-dimensional standard Gaussian vectors.
\State $x_i^j =l^m_{i}  + \sqrt{V^m_i }\  \Theta_i^j$
\State $x_{i+1}^j = x_i^j + b\left(t_i,x_i^j,\Yc^i_{\hat \theta_{h}^y}\left(x_i^j\right),\Zc^i_{\hat \theta_{h}^z}\left(x_i^j\right),\widetilde{u_{i}}\right) \Delta t + \sqrt{\Delta t}\ \sigma\left(t_{i},x_i^j,\widetilde{u_{i}}^{X}\right) \Xi_i^j$
\If{$i = N-1$}
\State $Y_{i+1}^j = g\left(x_N^j,\widetilde{u_{N}}^{X}\right)$
\Else\State $Y_{i+1}^j = \Yc^{i+1}_{ \theta_{m+1,i+1}}\left(x_{i+1}^j\right)$
\EndIf
\EndFor
\State $J^i(\hat \theta_{h}) = \frac{1}{B}\sum_{j=1}^{B} \bigg( f\left({t_{i}},x_i^j,\Yc^i_{\hat \theta_{h}^y}\left(x_i^j\right),\Zc^i_{\hat \theta_{h}^z}\left(x_i^j\right),\widetilde{u_{i}}\right)\ \Delta t +
Y_{i+1}^j - \Yc^i_{\hat \theta_{h}^y}\left(x_i^j\right)  -  \sqrt{\Delta t}\ \Zc^i_{\hat \theta_{h}^z}\left(x_i^j\right) \Xi_i^j \bigg)^2$
\State Calculate $\nabla J^i(\hat \theta_{h})$ by back-propagation.
\State Update $\hat \theta_{h+1} = \hat \theta_{h} - \rho_{h} \nabla J^i(\hat \theta_{h})$.
\EndFor
\State $\theta_{m+1,i}= \hat \theta_{H}$
\EndFor
\EndFor
\end{algorithmic}
\end{algorithm}

\section{Numerical results}\label{results}
The algorithms are implemented in Python with the Tensorflow library \cite{TF16}. Each numerical experiment is conducted using a node composed of 2 Intel® Xeon® Gold 5122 Processors, 192 Go of RAM, and 2 GPU nVidia® Tesla® V100 16Go. The multi-GPU parallelization on the global solver is conducted using the Horovod library \cite{horovod}.
The methods we test are:

\begin{itemize}
    \item \textbf{Global Direct}: algorithm \ref{algo: direct} at page \pageref{algo: direct}. Batch size $B = 10000$.
    \item \textbf{Global Dynamic}: algorithm \ref{algo: dynamic} at page \pageref{algo: dynamic}. Batch size $B = 200$ and $M = 100$.
    \item \textbf{Global Expectation}: algorithm \ref{algo: expectation} at page \pageref{algo: expectation}. Batch size $B = 2000 $.
    \item \textbf{Local Dynamic}: algorithm \ref{algo: local} at page \pageref{algo: local}. Batch size $B = 300\ (\mathrm{Weak}), B = 100\ (\mathrm{Pontryagin})$ and $M = 20 $, $R=50000$.
\end{itemize}
The name depends on the two choices made regarding the algorithm: the use a Global or Local method but also what law estimation technique is applied.
If the algorithm is applied to equations coming from the Pontryagin (abbreviated in Pont.) or the Weak approach, it is specified in its name. In the following tables, the incorrect results are highlighted in red. In Section \ref{sec: price impact} the three best results are highlighted in green whereas in Section \ref{sec: beyond} only the best value is.

\subsection{Linear price impact model}\label{sec: price impact}

We use a linear-quadratic mean-field game of controls model studied in \cite{AVLCDC19} and \cite{carmona2018probabilistic} for comparison. This model is useful for numerical tests since the analytic solution is known. The MFG of controls model for the representative player is given by:
\begin{equation}\label{eq: price impact}
\begin{aligned}
& \min_{\alpha\in\mathbb{A}}& &\E\left[\int_0^T \left(\frac{c_\alpha}{2}\norme{\alpha_t}^2 + \frac{c_X}{2} \norme{X_t}^2 - \gamma X_t \cdot u_t \right)\ \di t + \frac{c_g}{2} \norme{X_T}^2\ \right]\\
& \text{subject to} & &X_t = x_0 + \int_0^t \alpha_s\ \di s + \sigma\ W_t
\end{aligned}.
\end{equation} 
and the fixed point $\E[\alpha_t] =u_t$. In this case, the mean-field interaction is exerted through the law of the control process.\\
The Pontryagin optimality principle gives the system:
\begin{equation}\label{eq: weak price impact FBSDE}
\begin{cases}
\di X_t &= - \frac{1}{c_\alpha}Y_t\ \di t + \sigma\ \di W_t\\
X_0 &= x_0\\
\di Y_t &= - (c_X X_t + \frac{\gamma}{c_\alpha}\E[Y_t])\ \di t + Z_t\ \di W_t\\
Y_T &= c_g X_T.
\end{cases}
\end{equation}
In this case, the output $Z$ of the neural network is a matrix of size $d\times d$ and $Y$ is a vector of size $d$. \\
The weak representation of the value function gives:
\begin{equation}\label{eq: price impact FBSDE}
\begin{cases}
\di X_t &= - \frac{1}{c_\alpha}\sigma^{-1}Z_t\ \di t + \sigma\ \di W_t\\
X_0 &= x_0\\
\di Y_t &= - \left(\frac{c_X}{2} \norme{X_t}^2 + \frac{\gamma}{c_\alpha}X_t\cdot\sigma^{-1}\E[Z_t] + \frac{1}{2c_\alpha} \norme{\sigma^{-1} Z_t}^2\right)\ \di t + Z_t\ \di W_t\\
Y_T &= \frac{c_g}{2} \norme{X_T}^2.
\end{cases}
\end{equation}
 In this case, the output $Z$ of the neural network is a vector of size $d$ and $Y$ is a scalar. Therefore we may be able work in higher dimensions.
\begin{Rem}
With LQ models, the dynamics of $Y$ is linear in the Pontryagin approach and quadratic in the Weak approach. Thus the potentially high dimension of one method is counterbalanced by the complex dynamics of the other technique.
\end{Rem}

For our numerical experiments we take $c_X = 2, x_0 = 1, \sigma = 0.7, \gamma =2, c_\alpha = 2/3, c_g = 0.3$. If not stated otherwise, the simulations are conducted with $T = 1, d = 10, \Delta t = 0.01$. 

\begin{table}[H]
    \centering
    \begin{tabular}{|l||*{5}{c|}}\hline
    \backslashbox[25mm]{Method}{$T$}
    &\makebox[3em]{0.25}&\makebox[3em]{0.75}&\makebox[3em]{1.0}&\makebox[3em]{1.5}\\ \hline\hline
    \textbf{Reference}  & \textbf{0.7709}  & \textbf{0.1978} & \textbf{0.0811} & \textbf{0.0125}\\
    \hline
    Glob. Direct Pontryagin & 0.763  (1.3e-03) &
0.187  (2.5e-03) &
 0.075  (2.7e-03) &
\cellcolor{green} 0.012  (5.0e-03)\\ \hline
    Glob. Dyn. Pont. & 0.762  (2.3e-03) &
\cellcolor{green} 0.189  (4.0e-03) &
\cellcolor{green} 0.078  (5.5e-03) &
\cellcolor{green} 0.013  (6.7e-03) \\ \hline
Glob. Exp. Pont. ($0.1$) & 0.763 (1.6e-03) & \cellcolor{red}0.604 (1.1e-01) & \cellcolor{red} 0.729 (1.1e-01) & \cellcolor{red}0.803 (1.5e-01)\\ \hline
Glob. Exp. Pont. ($1.$) & 0.762 (1.4e-03) & \cellcolor{red} 0.251 (2.7e-02) & \cellcolor{red}0.467 (7.6e-02)  & \cellcolor{red}0.639 (1.1e-01)  \\ \hline
Glob. Exp. Pont. ($10.$) & 0.763 (1.5e-03) & 0.216  (1.7e-02) &
\cellcolor{red} 0.275  (3.7e-02) & \cellcolor{red}
0.574 (1.7e-01) \\ \hline
Glob. Exp. Pont. ($100.$) & \cellcolor{green}0.776 (8.4e-03) & \cellcolor{red}0.797 (1.1e-01) & \cellcolor{red}1.042 (1.3e-01)
 & \cellcolor{red}1.613 (2.6e-01)
 \\\hline
    Glob. Direct Weak & 0.778  (2.0e-03) &
\cellcolor{green}0.200  (1.4e-02) &
0.092  (2.9e-02) &
0.025  (2.0e-02)
\\ \hline
    Glob. Dyn. Weak & \cellcolor{green}0.775  (4.4e-03) &
0.212  (2.0e-02) &
\cellcolor{green} 0.083  (4.1e-02) &
0.016  (5.6e-02)
\\ \hline
Glob. Exp. Weak ($0.1$) & \cellcolor{red}0.877 (1.9e-02) & \cellcolor{red} 0.654 (9.9e-02) & \cellcolor{red} 0.595 (2.3e-01) & \cellcolor{red} 0.28 (6.0e-01)  \\ \hline
Glob. Exp. Weak ($1.$) &  \cellcolor{red} 0.901 (2.2e-03) & \cellcolor{red} 0.664 (9.8e-02) &  \cellcolor{red}  0.617 (1.0e-01) & \cellcolor{red} 0.507 (1.9e-01)  \\ \hline
Glob. Exp. Weak ($10.$) & \cellcolor{red} 0.887 (1.1e-02)   & \cellcolor{red} 0.698 (7.3e-02)  & \cellcolor{red}  0.6541 (6.3e-02)
  &  \cellcolor{red} 0.49 (2.3e-01) \\ \hline
Glob. Exp. Weak ($100.$) & \cellcolor{red} 0.887 (2.0e-03)  & \cellcolor{red} 0.650 (9.4e-02) &  \cellcolor{red} 0.602 (9.1e-02)
 &  \cellcolor{red} 0.492 (2.4e-01)
 \\ \hline
 Loc. Dyn. Pont & \cellcolor{green} 0.767 (3.5e-04) &
\cellcolor{green} 0.189 (6.3e-04) &
\cellcolor{green} 0.076 (7.6e-04) &
\cellcolor{green} 0.011 (7.5e-04)
 \\ \hline
 Loc. Dyn. Weak & \cellcolor{red} 0.944 (8.7e-04) &
\cellcolor{red}0.740 (2.6e-02) &
\cellcolor{red}0.692 (1.6e-02) &
\cellcolor{red}0.625 (2.2e-02)
 \\ \hline
    \end{tabular}
    \caption{Mean of $\E[X_T]$ over the 10 dimensions (and standard deviation) for several maturities $T$ on the price impact model \eqref{eq: price impact}. We perform 2000 iterations for global (glob.) methods and 20000 iterations for local (loc.) methods. For the expectation method, the value of the $\lambda$ penalization parameter is given under parenthesis.}
    \label{tab: expectations}
\end{table}

\begin{table}[H]
    \centering
    \begin{tabular}{|l||*{2}{c|}}\hline
    Global Direct Pontryagin & \textbf{1877 s.} \\\hline
    Global Dyn. Pontryagin & \textbf{1336 s.} \\\hline
    Global Exp. Pontryagin & \textbf{1562 s.}
 \\\hline
    Global Direct Weak & \textbf{2205s.}    \\\hline
    Global Dyn. Weak &  \textbf{1605 s.} \\\hline
    Global Exp. Weak &  \textbf{1670 s.} \\\hline
    Local Dynamic Pontryagin & 11627 s.
 \\\hline
    Local Dynamic Weak &  12689 s.\\\hline
    \end{tabular}
    \caption{Duration times of the methods (2000 iterations for global methods, 20000 iterations for local methods) on the price impact model \eqref{eq: price impact} with $T=1$. on one run}
    \label{tab: duration trader}
\end{table}

\begin{figure}[H]
   \begin{minipage}[c]{.49\linewidth}
      \includegraphics[width=0.9\linewidth]{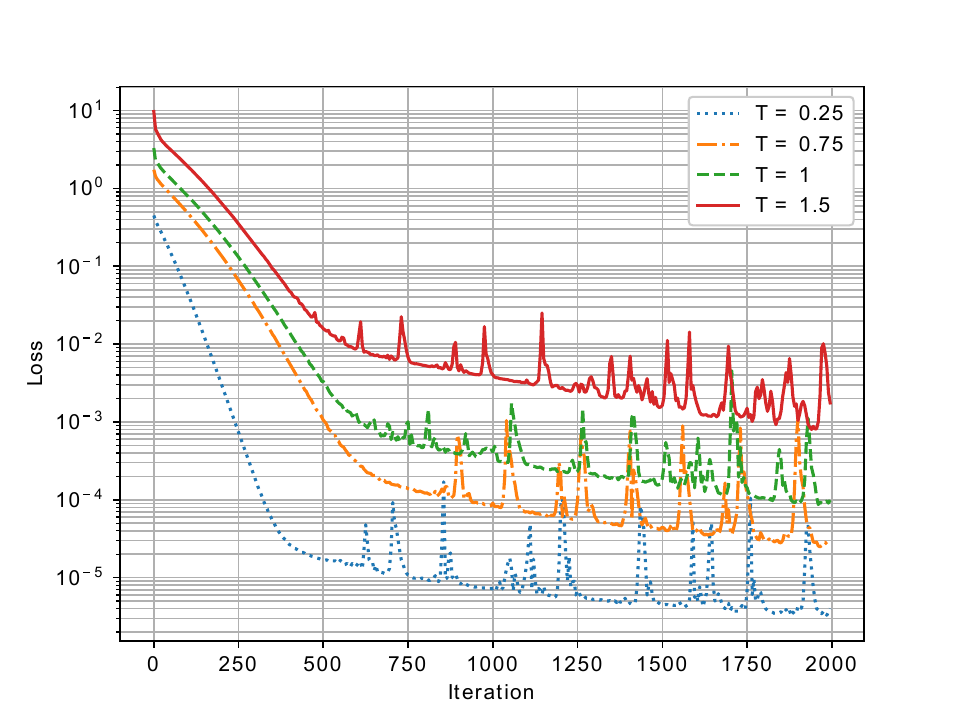}
   \end{minipage} \hfill
   \begin{minipage}[c]{.49\linewidth}
      \includegraphics[width=0.9\linewidth]{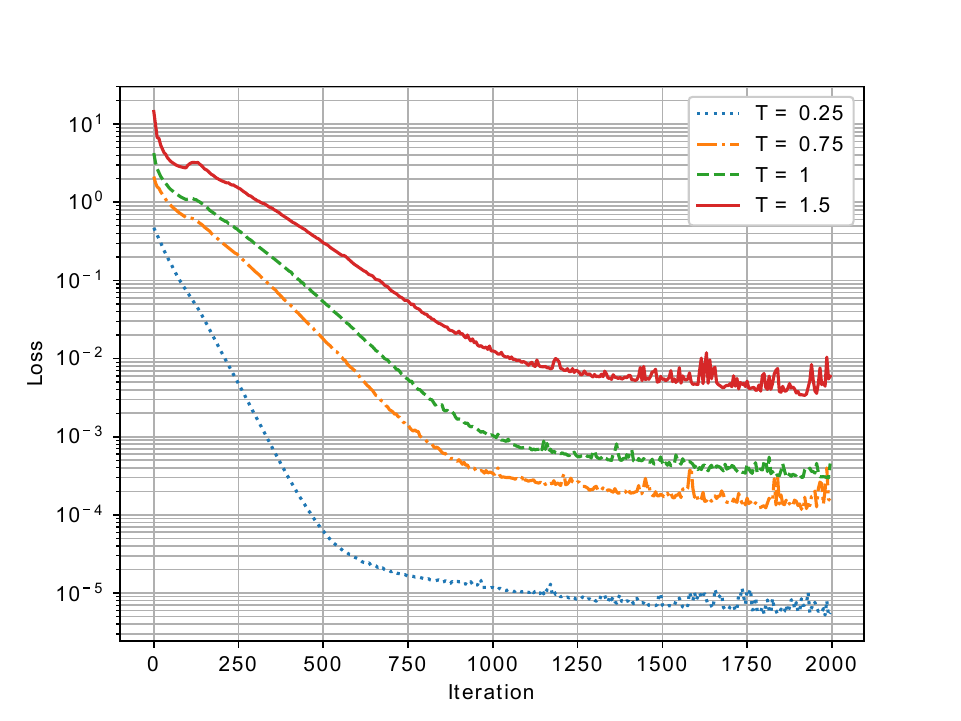}
   \end{minipage}
    \caption{Learning curves for Global Direct (left) and Global Dynamic (right) Pontryagin methods on the price impact model \eqref{eq: price impact}. The loss is the $L^2$ error between $Y_T$ and the terminal condition of the backward equation.}
   \label{fig: learning Pont}
\end{figure}

\begin{figure}[H]
   \begin{minipage}[c]{.49\linewidth}
      \includegraphics[width=0.9\linewidth]{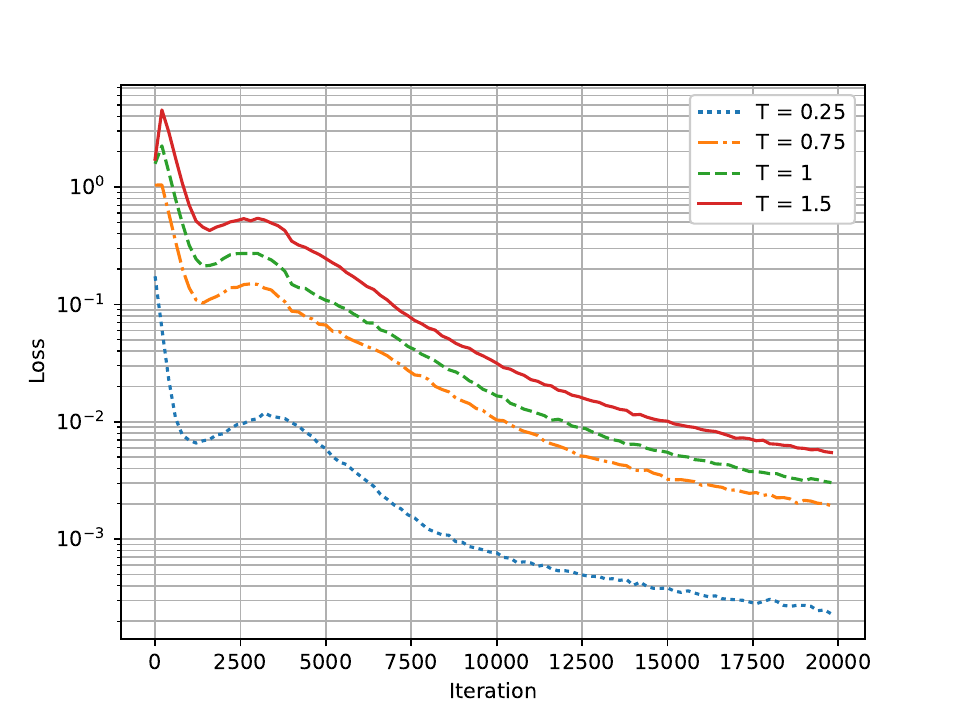}
   \end{minipage} \hfill
   \begin{minipage}[c]{.49\linewidth}
      \includegraphics[width=0.9\linewidth]{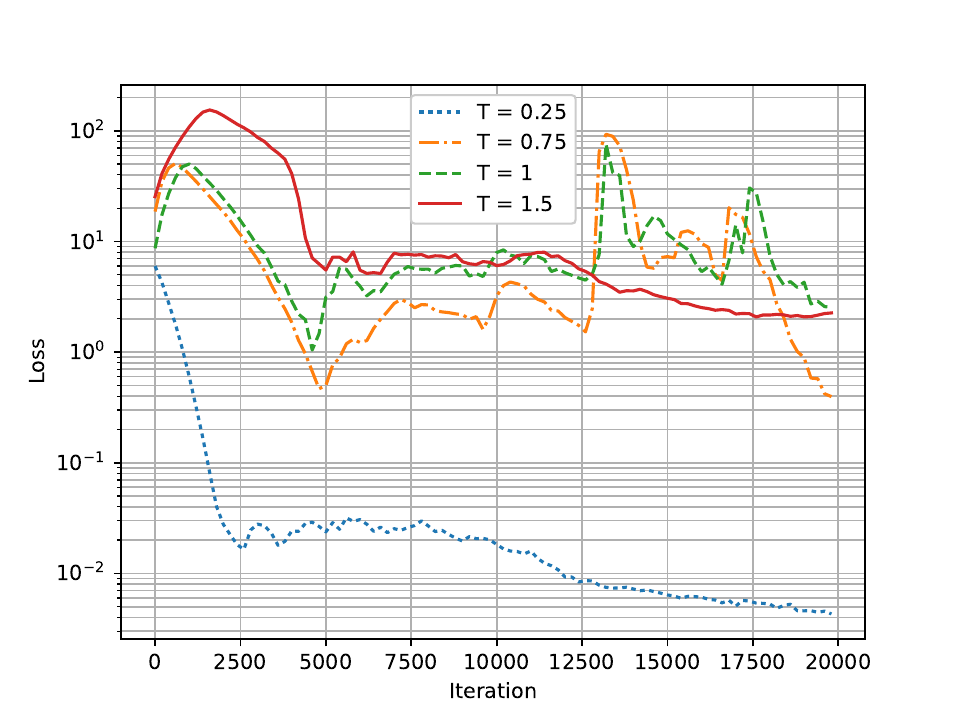}
   \end{minipage}
   \caption{Learning curves for Local Dynamic Pontryagin (left) and Local Dynamic Weak (right) methods on the price impact model \eqref{eq: price impact}. The loss is the sum of the local $L^2$ errors between the neural network $Y$ and the Euler discretization for all time steps.}
    \label{fig: learning curve local}
\end{figure}

\begin{figure}[H]
   \begin{minipage}[c]{.49\linewidth}
      \includegraphics[width=0.9\linewidth]{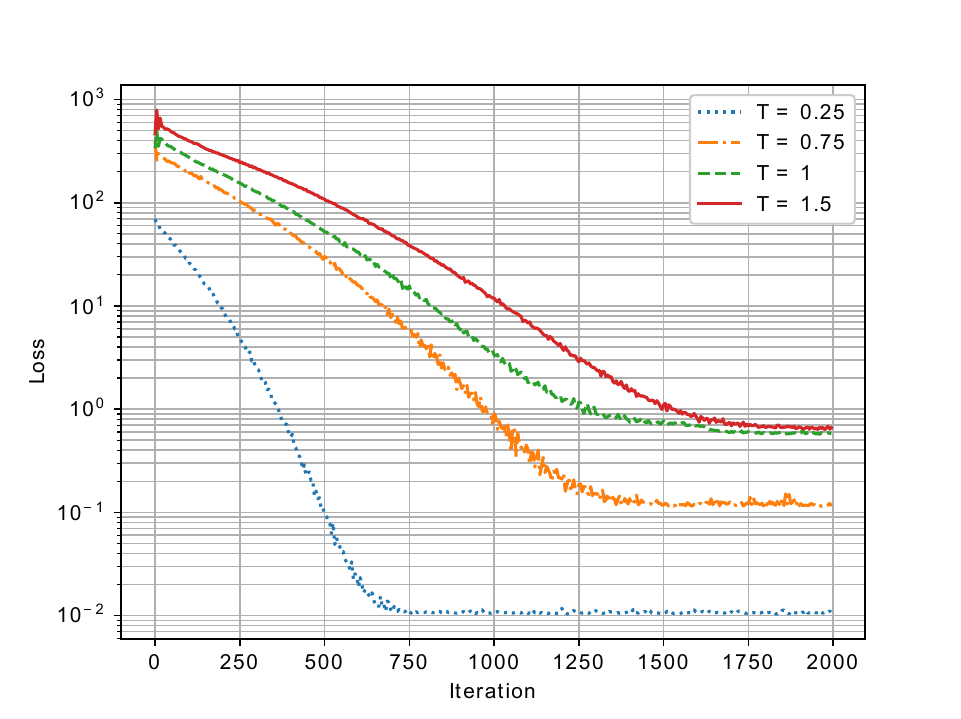}
   \end{minipage} \hfill
   \begin{minipage}[c]{.49\linewidth}
      \includegraphics[width=0.9\linewidth]{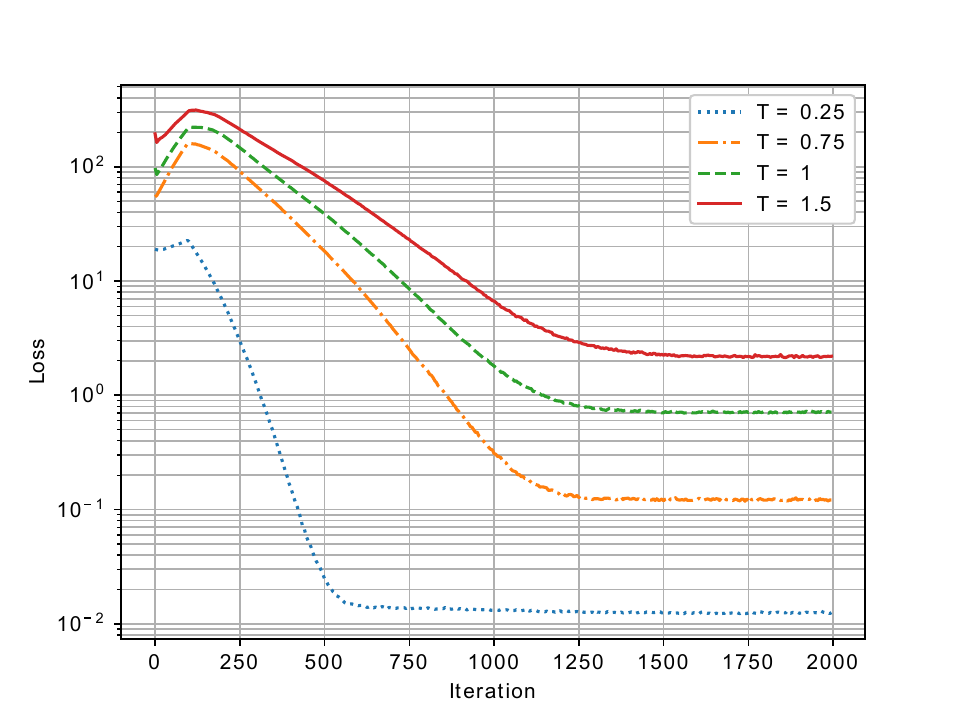}
   \end{minipage}
    \caption{Learning curves for Global Direct Weak (left) and Global Dynamic Weak (right) methods on the price impact model \eqref{eq: price impact}. The loss is the $L^2$ error between $Y_T$ and the terminal condition of the backward equation.}
   \label{fig: learning curve weak}
\end{figure}

\begin{figure}[H]
   \begin{minipage}[c]{.49\linewidth}
      \includegraphics[width=0.9\linewidth]{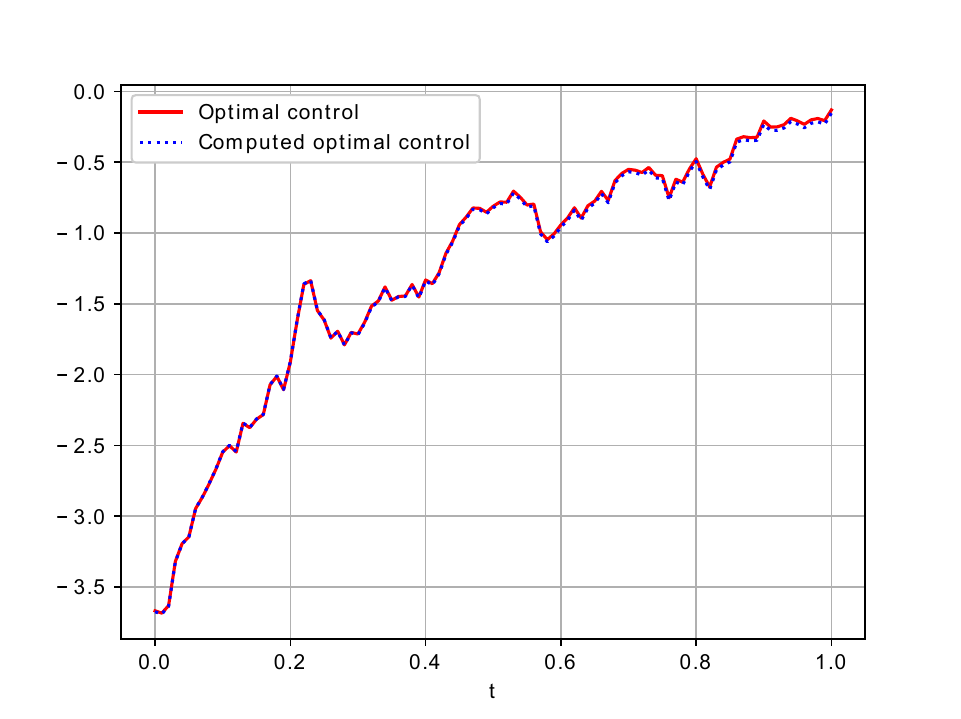}
   \end{minipage} \hfill
   \begin{minipage}[c]{.49\linewidth}
      \includegraphics[width=0.9\linewidth]{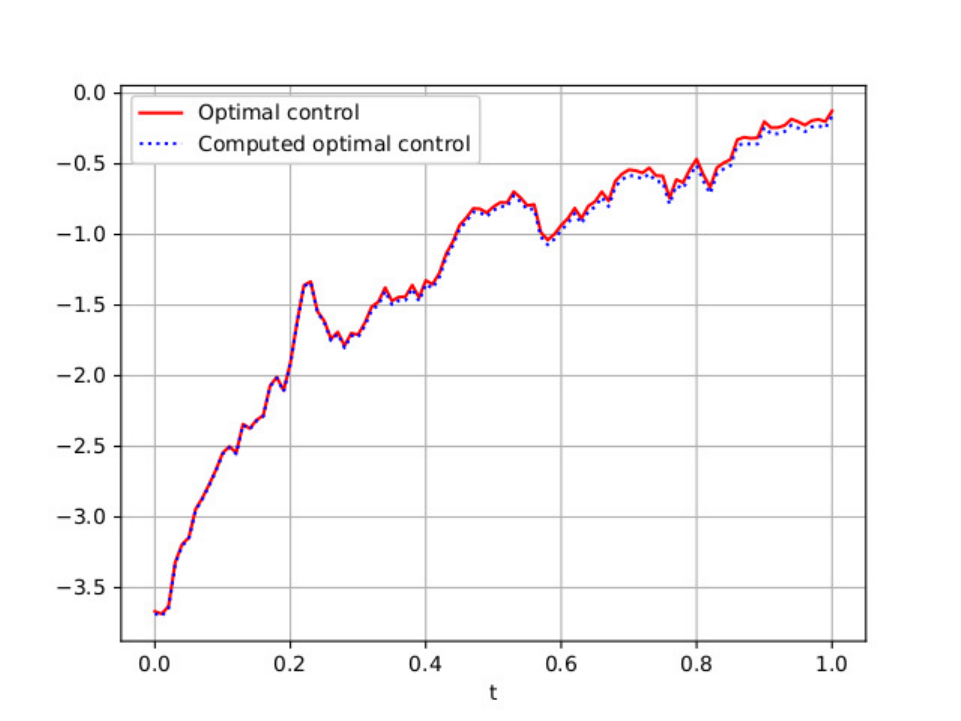}
   \end{minipage}
   \caption{First coordinate of the optimal control evaluated on a sample path for Global Direct Pontryagin (left) and Global Dynamic (right) Pontryagin methods after 2000 iterations on the price impact model \eqref{eq: price impact} with $T=1$.}
   \label{fig: control Pont}
\end{figure}

\begin{figure}[H]
   \begin{minipage}[c]{.49\linewidth}
      \includegraphics[width=0.9\linewidth]{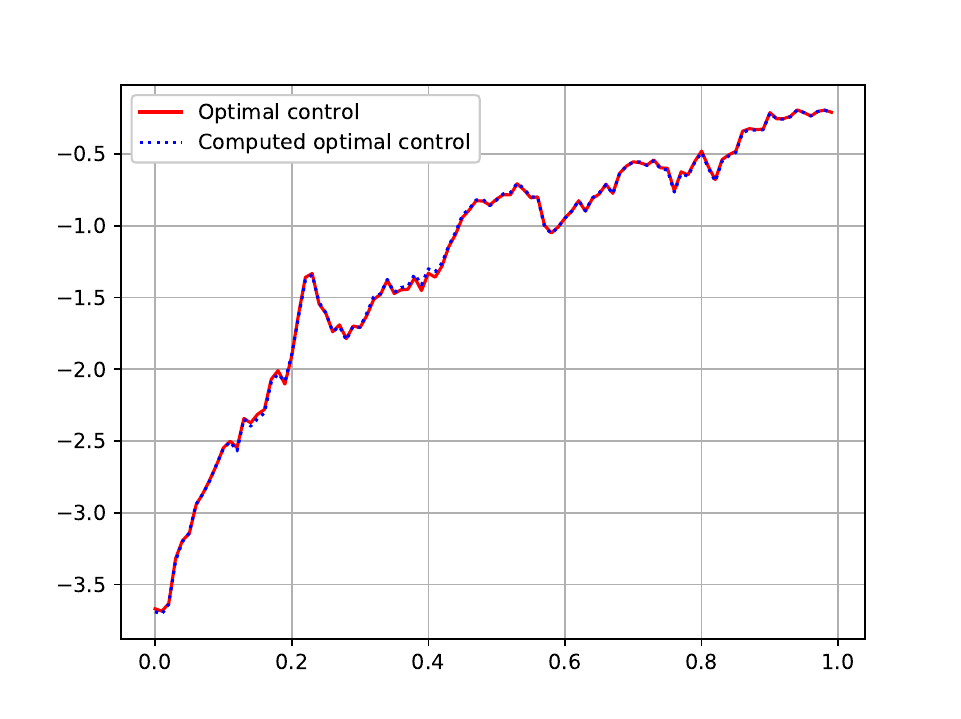}
   \end{minipage} \hfill
   \begin{minipage}[c]{.49\linewidth}
      \includegraphics[width=0.9\linewidth]{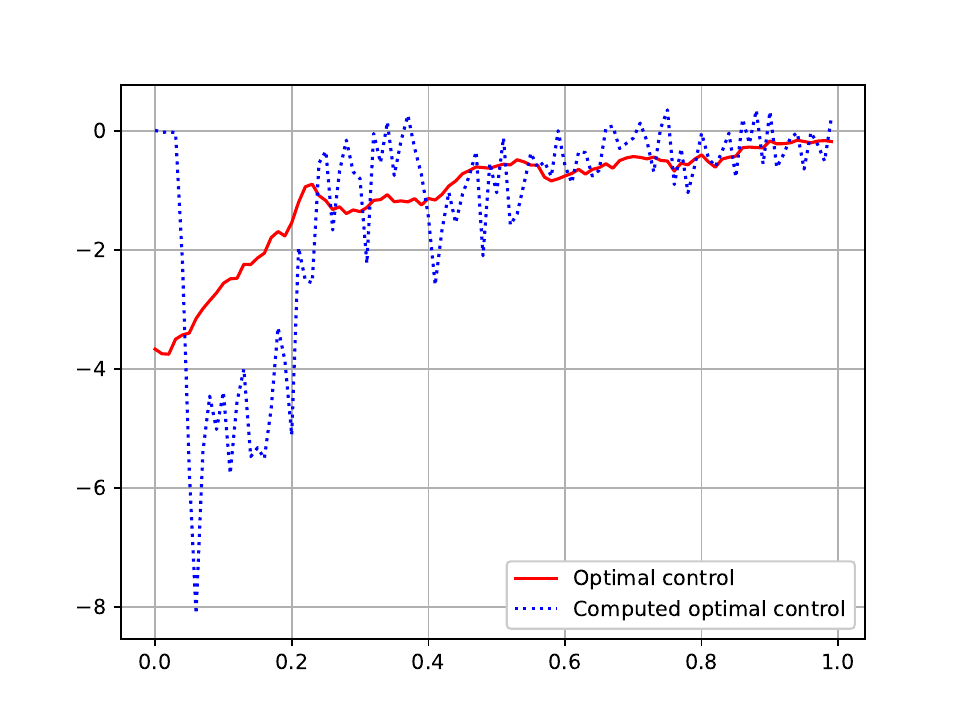}
   \end{minipage}
   \caption{First coordinate of the optimal control evaluated on a sample path for Local Dynamic Pontryagin (left) and Local Dynamic Weak (right) methods after 20000 iterations on the price impact model \eqref{eq: price impact} with $T=1.$}
   \label{fig: control local}
\end{figure}

\paragraph{}
All methods except Global Expectation Weak and Local Dynamic Weak converge to the exact solution for small maturities.
These two solvers do not converge to the right solution for any time horizon.

The choice of the parameter $\lambda$ influences a lot the output of the Global Expectation Pontryagin scheme, as observed in Table \ref{tab: expectations}. The best results are obtained when $\lambda$ is of order $10$ but we notice that a large range of values seems to work fine for very small time horizon $T$. However the Global Expectation Weak scheme never works. 
We will not test the expectation methods on the other test cases since they are less efficient than the other methods.

We see in Figure \ref{fig: learning curve local} that the Local Dynamic Pontryagin  method needs more iterations for the loss to stabilize than the Global methods. We cannot hope for more iterations to help the convergence in the Global Weak methods since the loss in the learning curves of Figure \ref{fig: learning curve weak} reaches a plateau. The algorithms solving the system coming from the Pontryagin principle perform better than the others. The dynamic estimation of the expectation allows to gain training speed and to smooth the loss, as seen in Figure \ref{fig: learning Pont} and Table \ref{tab: duration trader}. As another accuracy test, we can also plot the optimal control for which we have an analytical expression. We see in Figures \ref{fig: control Pont}, \ref{fig: control local} that Global Pontryagin and Local Dynamic Pontryagin methods perform well but that the Local Dynamic Weak method does not seem to converge, which confirms what is observed in Table \ref{tab: expectations}.

\subsection{A one-dimensional mixed model}\label{pop}
\begin{figure}[H]
    \centering
    \includegraphics[width=10cm]{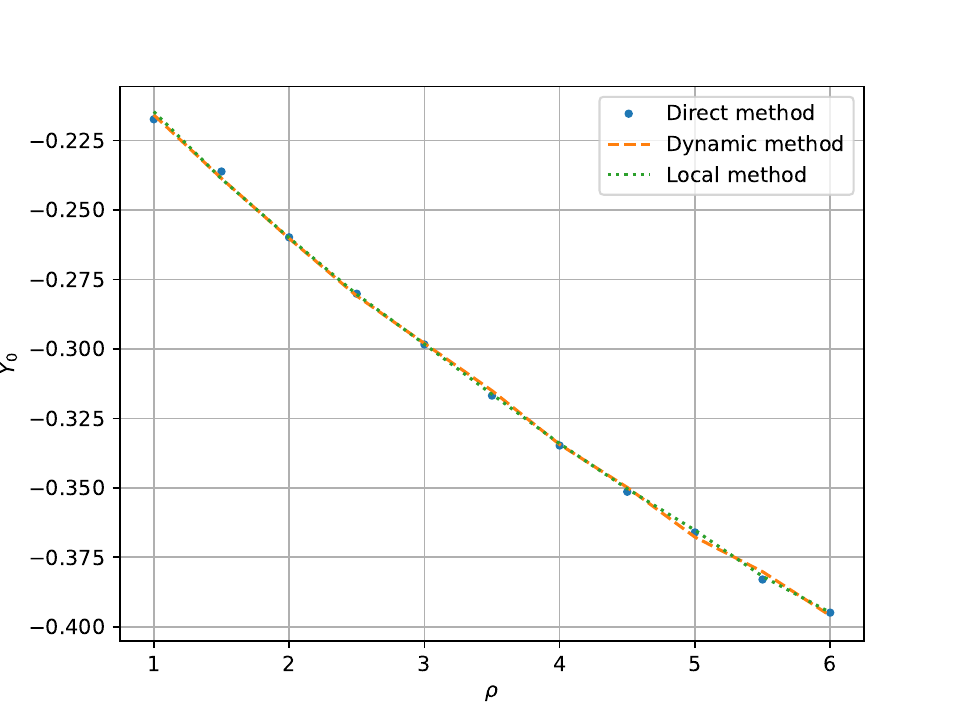}
    \caption{Value of $Y_0$ as a function of parameter $\rho$ for the model \eqref{eq: population}}
    \label{fig: pop}
\end{figure}

We consider the following one-dimensional example from \cite{CL19b,AVLCDC19}:
\begin{align}
\begin{cases}
    \di X_t & = - \rho Y_t\ \di t + \sigma\ \di W_t,\ X_0 = x_0\\
    \di Y_t & = \arctan(\E[X_t])\ \di t + Z_t\ \di W_t,\ Y_T = \arctan(X_T).\label{eq: population}
\end{cases}
\end{align} This model comes from the Pontryagin principle applied to the mean-field game problem
\begin{align*}
    \min_\alpha\ &\E\Big[\int_0^T \big(\frac{1}{2\rho} \alpha_s^2 - X_s \arctan( u_s)\big)\ \di s + g(X_T)\Big]\\
     & \di X_t  = \alpha_t\ \di t + \di W_t,\ X_0 = x_0,
\end{align*} with the fixed point $u_s = \E[X_s]$, and where $g$ is an antiderivative of $\arctan$. 
We take the same model parameters as in \cite{CL19b} ($T=1$ and $x_0=1$) and obtain in Figure \ref{fig: pop} with all our methods the same results as in their Figure 4. For the numerical resolution we choose $100$ time steps. Notice that we use 3 hidden layers with 11 neurons in each when \cite{CL19b} uses 100 neurons by layer. We see that our smaller number of neurons is enough for this example resolution.
The duration of the methods are given in Table \ref{tab: pop} which illustrates again the speed gain in using the dynamic method.
\begin{table}[t]
    \centering
    \begin{tabular}{|l||*{5}{c|}}
    \hline
        Global Direct & \textbf{535 s.}\\
        \hline
        Global Dynamic & \textbf{425 s.}\\
        \hline
        Local Dynamic & \textbf{10900 s.}\\
        \hline
    \end{tabular}
    \caption{Duration times of the algorithms for model \eqref{eq: population} on one run (2000 iterations for global methods, 20000 iterations for the local method)}
    \label{tab: pop}
\end{table}

\subsection{Beyond the mean-field games case}\label{sec: beyond}
\paragraph{}
In this Section we design non Linear Quadratic models in order to test the limitations of our methods. We construct general MKV FBSDES with explicit solutions following a log-normal distribution. Let $X_t^i$ be defined by

\begin{align}
\label{forward}
  \di X_t^i & = a X^i_t\ \di t + \sigma X^i_t\ \di W^i_t, \\
  X^i_0 & = \xi,
\end{align} with $a$, $\sigma$, $\xi \in \R$. We obtain explicitly
\begin{align*}
  X^i_t & = \xi e^{(a-\frac{\sigma^2}{2}) t +\sigma W^i_t}, \\
   h^i_t & :=  \E[ X^i_t] = \xi e^{a t},  \\
    k^i_t & := \E[ (X^i_t)^2]  = \xi e^{(2a +\sigma^2) t }. 
\end{align*}

We choose  $u: (t,x) \mapsto e^{\alpha t} \log\left( \prod_{i=1}^n x^i \right)$  and the following dynamic for $Y_t$
\begin{align*}
  Y_t & = e^{\alpha t} \log\left(\prod_{i} X^i_t\right) = e^{\alpha t}  \sum_i \left[ \log(\xi)+ (a-\frac{\sigma^2}{2}) t + \sigma W^i_t \right],
  \end{align*}
  such that
  \begin{align*}
  c_t & :=\E[ Y_t ]  = e^{\alpha t}  \sum_i \left[\log(\xi)+ \left(a-\frac{\sigma^2}{2}\right) t\right],\\
  d_t & := \E[Y_t^2]  = e^{2 \alpha t} \left( \left[\sum_i  \left( \log(\xi)+ a-\frac{\sigma^2}{2}\right) t\right]^2 + \sum_i \sigma^2 t \right).
\end{align*}
As we have $Y_t = u(t,X_t)$, we obtain $Z_t^i = \sigma X^i_t \partial_x u(t,X_t)$ following:
\begin{align*}
  Z_t^i & = \sigma e^{\alpha t} \\
  e_t^i & :=\E[Z_t^i ]  = \sigma e^{\alpha t},\\
  f_t^i & := \E[(Z_t^i)^2 ]  = \sigma^2 e^{2\alpha t}.
\end{align*}
Introducing 
\begin{align*}
 \phi(t,x) & := \partial_t u + \sum_{i} a x_i\ \partial_{x_i} u + 
 \sum_{i} \frac{(\sigma x_i)^2}{2}\ \partial_{x_i^2}^2 u  \\
  & = e^{\alpha t} \left( \alpha \log\left(\prod_{i} x^i\right)  + \sum_i \left(a - \frac{\sigma^2}{2}\right)\right),
\end{align*}
$u(t,X_t)$ solves the PDE
\begin{align*}
\partial_t u + \sum_{i} a x^i\ \partial_{x_i} u + 
\sum_{i} \frac{\sigma^2}{2}\ \partial_{x_i^2}^2 u - \phi(t,x)= 0.
\end{align*}
This semilinear PDE is related to the BSDE associated with the driver $f(t,x) = -\phi(t,x)$ for forward dynamics \eqref{forward}.\\
Using some chosen $\R^d$ valued functions $\psi$ and $\R^k$ valued functions $\kappa$,
 we express all dynamics in a McKean-Vlasov setting:
\begin{equation}
\begin{cases}
\di X_t^i  & = (a X^i_t +  \psi(Y_t, Z_t^i, \E[X_t^i],\E[(X_t^i)^2], \E[Y_t], \E[Y_t^2],
\E[Z_t^i], \E[(Z^i_t)^2])\\ & - \psi\left(e^{\alpha t} \log\left(\prod_{i} X^i_t\right), \sigma e^{\alpha t}, h_t^i,k_t^i,c_t,d_t,e_t^i,f_t^i\right)\ \di t + \sigma X^i_t\ \di W^i_t  \\
  X^i_0 & = \xi\\
\di Y_t & = - f(t,X_t,Y_t,Z_t,\E[X_t],\E[X_t^2], \E[Y_t], \E[Y_t^2],
\E[Z_t], \E[Z^2_t])\ \di t \\ & + Z_t\ \di W_t\\
Y_T & = e^{\alpha T} \log\left(\prod_{i} X^i_T\right)
\end{cases}    
\end{equation}
with
\begin{align*}
& f(t,X_t,Y_t,Z_t,x_1,x_2,y_1,y_2,z_1,z_2 ) \\& = -\phi(t,x) + \kappa(Y_t, Z_t,x_1,x_2,y_1,y_2,z_1,z_2 ) \\ &-
\kappa\left(e^{\alpha t} \log\left(\prod_{i} X^i_t\right), \sigma e^{\alpha t}, h_t^i,k_t^i,c_t,d_t,e_t^i,f_t^i\right).
\end{align*}
and $f: \R\times\R^d\times\R\times\R^{1\times d}\times\R^d\times\R^d\times\R\times\R\times\R^{1\times d}\times\R^{1\times d}\mapsto \R$.\\

We consider two models of this kind for numerical tests.

\subsubsection{A fully coupled linear example }
We consider a linear McKean-Vlasov FBSDE  in $Y_t$, $Z_t$ and their law dynamics for $X_t$ and $Y_t$:
\begin{equation}
\begin{cases}
\di X_t^i  & = (a X^i_t +  b(Y_t + Z_t^i + \E[X_t^i] + \E[Y_t] +
\E[Z_t^i]) \\ &  - b\left(e^{\alpha t} \log\left(\prod_{i=1}^d X^i_t\right)+ \sigma e^{\alpha t} + h_t^i + c_t + e_t^i\right)\ \di t  + \sigma X^i_t\ \di W^i_t  \\
  X^i_0 & = \xi\\
\di  Y_t & =  \bigg(\phi(t,X_t) + b(Y_t + \frac{1}{d}\sum_{i=1}^d Z_t^i + \frac{1}{d}\sum_{i=1}^d\E[X_t^i] + \E[Y_t] \\ &+ \frac{1}{d}\sum_{i=1}^d \E[Z_t^i])   - b\left(e^{\alpha t} \log\left(\prod_{i=1}^d X^i_t\right) + \frac{1}{d}\sum_{i=1}^d \sigma e^{\alpha t} + \frac{1}{d}\sum_{i=1}^d h_t^i \right. \\ &+ \left. c_t +  \frac{1}{d}\sum_{i=1}^d e_t^i\right) \bigg)\ \di t + Z_t\ \di W_t\\
Y_T & = e^{\alpha T} \log\left(\prod_{i=1}^d X^i_T\right).\label{eq: linear model}
\end{cases}    
\end{equation}

\begin{table}[t]
    \centering
    \begin{tabular}{|l||*{5}{c|}}\hline
    \backslashbox[25mm]{Method}{$T$}
    &\makebox[3em]{0.25}&\makebox[3em]{0.75}&\makebox[3em]{1.0}&\makebox[3em]{1.5}\\ \hline\hline
    \textbf{Reference}  & \textbf{1.0253}  & \textbf{1.0779} & \textbf{1.1052} & \textbf{1.1618}\\
    \hline
  Global  Direct & \cellcolor{green}\cellcolor{green}1.025 (2e-03) &
1.076 (3e-03) &
1.095 (4e-03) &
\cellcolor{green}1.162 (7e-03)
\\ \hline
 Global   Dynamic & 1.026 (2e-03) &
\cellcolor{green}1.077 (4e-03) &
\cellcolor{green}1.105 (3e-03) &
1.163 (5e-03)
 \\ \hline   Local Dynamic & 1.025 (2e-04)&
1.092 (5e-04)&
1.146 (8e-04) &
\cellcolor{red}1.28 (1e-03)
 \\ \hline
    \end{tabular}
    \caption{Mean of $\E[X_T]$ over the 10 dimensions (and standard deviation) for several maturities $T$ (2000 iterations for global methods, 20000 iterations for local method) on the fully coupled linear model \eqref{eq: linear model}.}
    \label{tab: res linear }
\end{table}

\begin{table}[t]
    \centering
    \begin{tabular}{|l|*{2}{c|}}\hline
    Global Direct & Global Dynamic & Local Dynamic\\ \hline
    \textbf{2081 s.} &
    \textbf{1308 s.}  &
    \textbf{14811 s.}    \\ \hline
    \end{tabular}
    \caption{Duration times of the methods (2000 iterations for global methods, 20000 iterations for local method) on the fully coupled linear model \eqref{eq: linear model} for $T=1$.}
\end{table}

\begin{figure}[H]
    \centering
   \includegraphics[width=.49\linewidth]{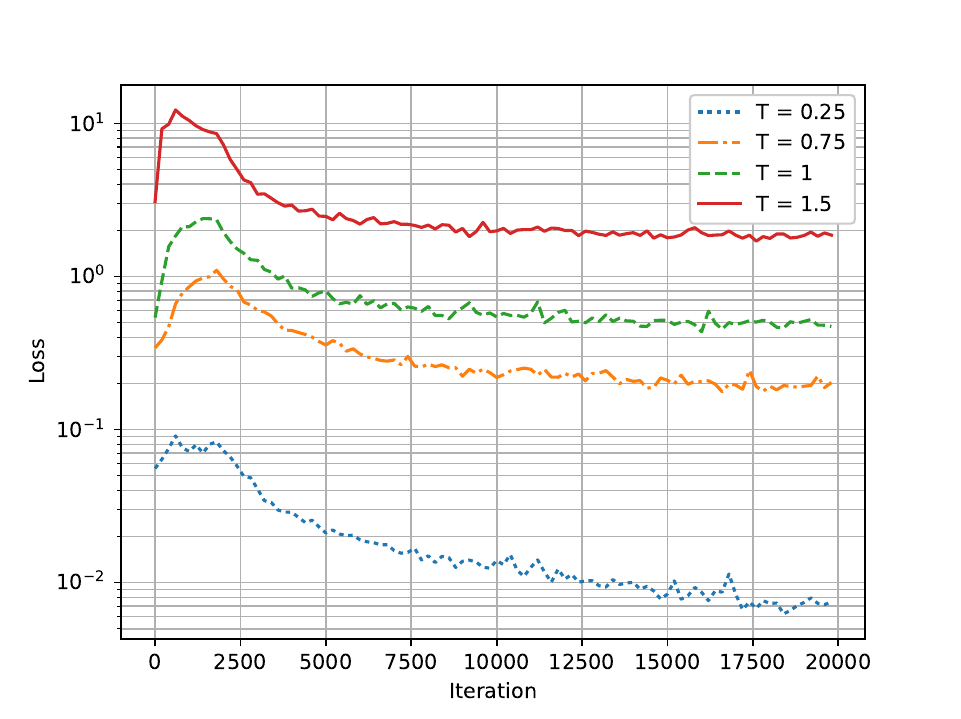}
   \caption{Learning curves for Local Dynamic method on the fully coupled linear model \eqref{eq: linear model}. The loss is the sum of the local $L^2$ errors between the neural network $Y$ and the Euler discretization for all time steps.}
   \label{fig: learning curve lin loc}
\end{figure}

\begin{figure}[H]
   \begin{minipage}[c]{.49\linewidth}
      \includegraphics[width=0.9\linewidth]{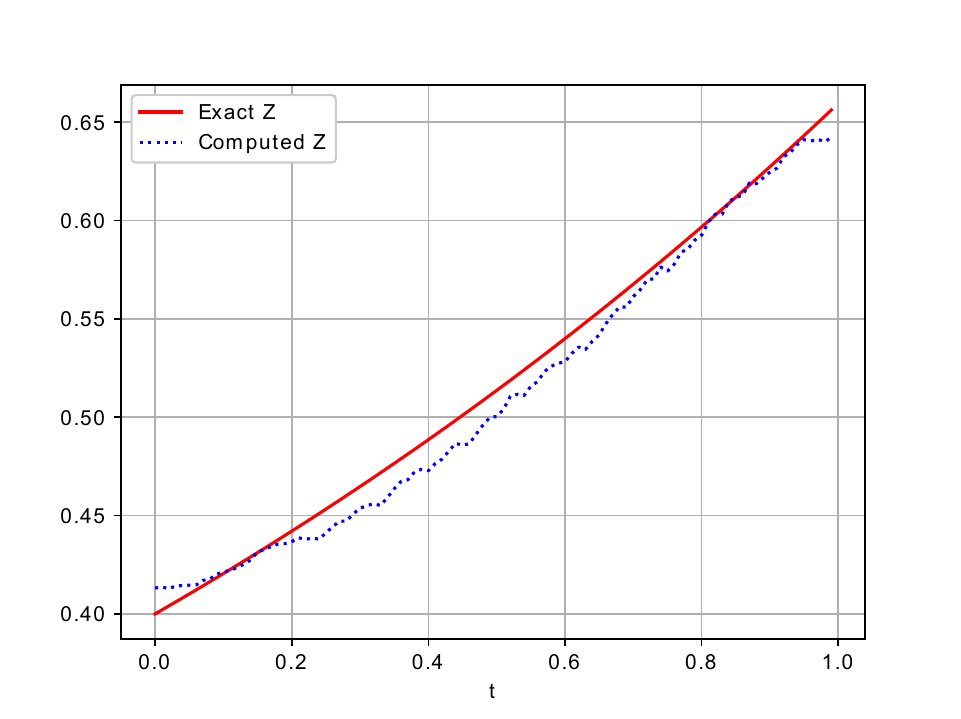}
   \end{minipage} \hfill
   \begin{minipage}[c]{.49\linewidth}
      \includegraphics[width=0.9\linewidth]{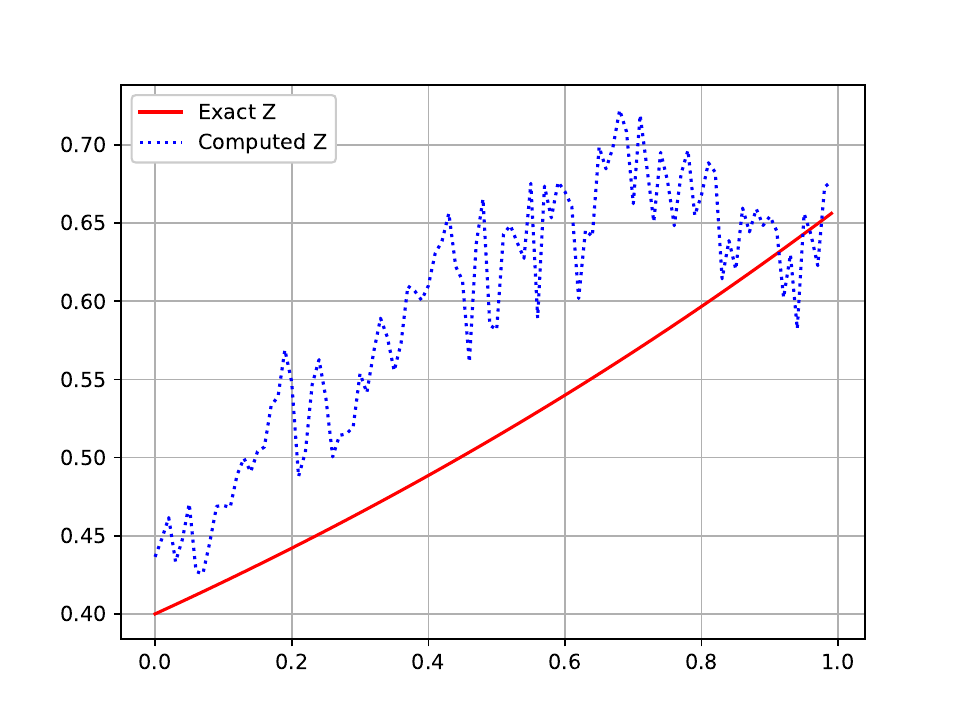}
   \end{minipage}
   \caption{First coordinate of $Z_t$ evaluated on a sample path for 
   Global Dynamic  (left) and 
   Local Dynamic (right) methods after 2000 iterations (respectively 20000) on the fully coupled linear model \eqref{eq: linear model}.}
   \label{fig: z lin loc}
\end{figure}

\begin{figure}[H]
   \begin{minipage}[c]{.49\linewidth}
      \includegraphics[width=0.9\linewidth]{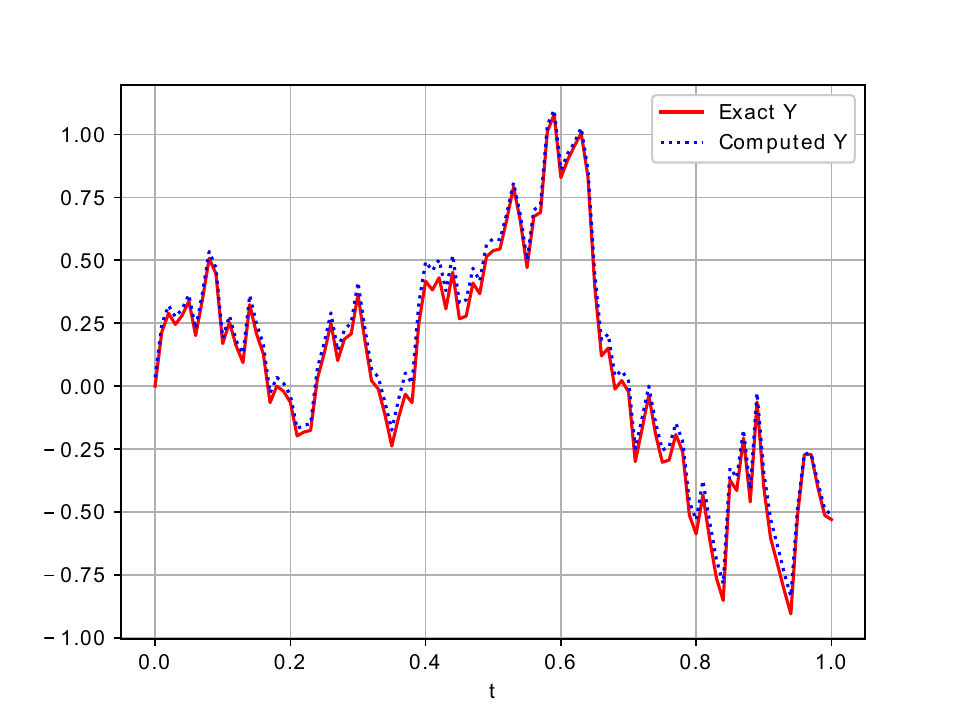}
   \end{minipage} \hfill
   \begin{minipage}[c]{.49\linewidth}
      \includegraphics[width=0.9\linewidth]{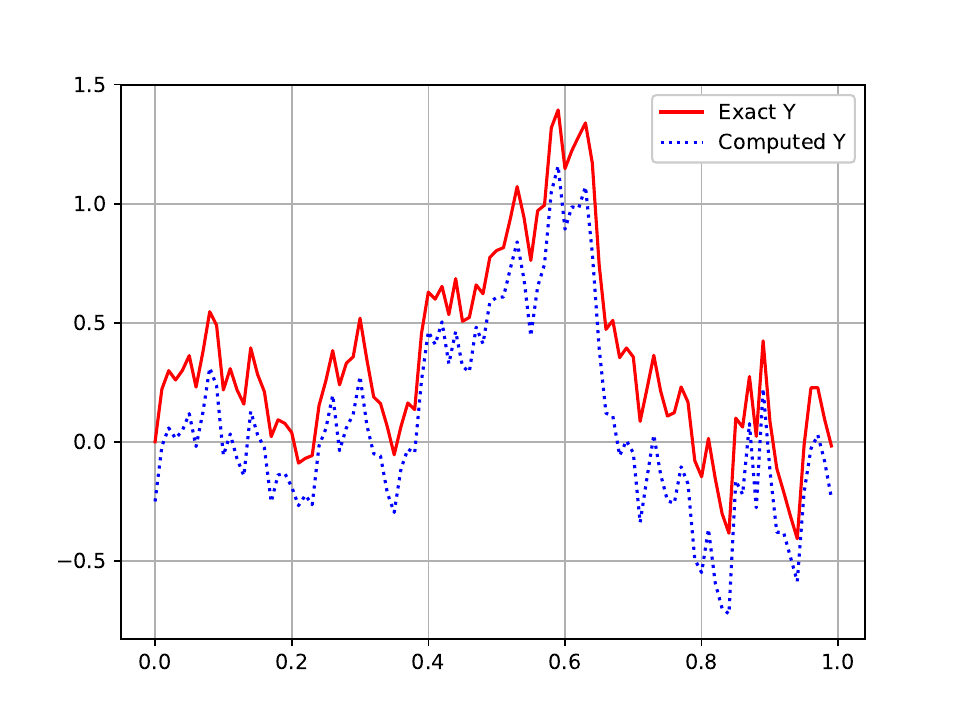}
   \end{minipage}
   \caption{First coordinate of $Y_t$ evaluated on a sample path for Global Direct (left) and Local Dynamic methods (right) 
   after 2000 iterations (respectively 20000) on the fully coupled linear model \eqref{eq: linear model}. }
   \label{fig: y lin loc}
\end{figure}

\paragraph{}
We take $a=b=0.1,\ \alpha = 0.5,\ \sigma = 0.4,\ \xi = 1$. The three algorithms demonstrate good performances on this test case. Both processes $Y, Z$ are well represented by the neural network. However the Local Dynamic method is less precise than the global methods when the maturity grows. We see in Table \ref{tab: res linear }, Figure \ref{fig: z lin loc}, Figure \ref{fig: y lin loc} that the Local Dynamic method is biased when $T=1$ when the Global methods achieve a great accuracy. It looks like the results of the Local Dynamic method cannot be improved since the loss flattens in Figure \ref{fig: learning curve lin loc}.

\subsubsection{A fully coupled quadratic example }
We consider a quadratic McKean-Vlasov FBSDE in $Y_t$, $Z_t$ and their law dynamics for $X_t$ and $Y_t$:
We consider a quadratic McKean-Vlasov FBSDE in $Y_t$, $Z_t$ and their law dynamics for $X_t$ and $Y_t$:
\begin{equation}    
\begin{cases}
\di X_t^i  & = (a X^i_t +  b(Y_t + Z_t^i + \E[X_t^i] + \E[Y_t] +
\E[Z_t^i]) + \sigma X^i_t\ \di W^i_t  \\ &  - b\left(e^{\alpha t} \log\left(\prod_{i=1}^d X^i_t\right)+ \sigma e^{\alpha t} + h_t^i + c_t + e_t^i\right) \\ &  +  c \Bigg(Y_t^2 + (Z_t^i)^2 + \E[(X_t^i)^2] + \E[Y_t^2] +
\E[(Z_t^i)^2]) \\ &  - c\left(e^{2\alpha t} \log\left(\prod_{i=1}^d X^i_t\right)^2+ \sigma^2 e^{2\alpha t} + (h_t^i)^2 + c_t^2 + (e_t^i)^2\right) \Bigg)\ \di t   \\
  X^i_0 & = \xi\\
\di Y_t & = \Bigg(\phi(t,X_t) + b(Y_t + \frac{1}{d}\sum_{i=1}^d Z_t^i + \frac{1}{d}\sum_{i=1}^d\E[X_t^i] + \E[Y_t]) \\ & + \frac{b}{d}\sum_{i=1}^d \E[Z_t^i])   - b\left(e^{\alpha t} \log\left(\prod_{i=1}^d X^i_t\right) + \frac{1}{d}\sum_{i=1}^d \sigma e^{\alpha t}\right)  + c \E[Y_t^2]\\ & -b\left( \frac{1}{d}\sum_{i=1}^d h_t^i + c_t +  \frac{1}{d}\sum_{i=1}^d e_t^i\right) + c(Y_t^2 + \frac{1}{d}\sum_{i=1}^d (Z_t^i)^2 + \frac{1}{d}\sum_{i=1}^d\E[(X_t^i)^2] ) \\ &+ \frac{c}{d}\sum_{i=1}^d \E[(Z_t^i)^2])  - c\left(e^{2\alpha t} \log\left(\prod_{i=1}^d X^i_t\right)^2 + \frac{1}{d}\sum_{i=1}^d \sigma^2 e^{2\alpha t}\right) \\ &- c\left( \frac{1}{d}\sum_{i=1}^d (h_t^i)^2 + c_t^2+ \frac{1}{d}\sum_{i=1}^d (e_t^i)^2\right)\Bigg)\ \di t + Z_t\ \di W_t\\
Y_T & = e^{\alpha T} \log\left(\prod_{i=1}^d X^i_T\right).\label{eq: quadratic model}
\end{cases}    
\end{equation}  

\begin{table}[t]
    \centering
    \begin{tabular}{|l||*{5}{c|}}\hline
    \backslashbox[25mm]{Method}{$T$}
    &\makebox[3em]{0.25}&\makebox[3em]{0.75}&\makebox[3em]{1.0}&\makebox[3em]{1.5}\\ \hline\hline
    \textbf{Reference}  & \textbf{1.0253}  & \textbf{1.0779} & \textbf{1.1052} & \textbf{1.1618}\\
    \hline
    Global Direct & 1.024  (1.8e-03) &
1.065  (4.3e-03) &
\cellcolor{red}12.776 (3.3e-02) & \cellcolor{red}DV
\\ \hline
   Global Dynamic & \cellcolor{green}1.025  (2.1e-03) &
\cellcolor{green}1.072 (3.1e-03) &
0.961  (7.0e-03) & \cellcolor{red}DV
 \\ \hline
    Local Dynamic & 1.024 (1.6e-04) &
\cellcolor{red}-7.180 (9.0e-04) &
0.411 (1.1e-03) 
 & \cellcolor{red}DV
 \\ \hline
    \end{tabular}
    \caption{Mean of $\E[X_T]$ over the 10 dimensions (and standard deviation) for several maturities $T$ (2000 iterations for global methods, 20000 iterations for Local method) on the fully coupled quadratic model \eqref{eq: quadratic model}.}
    \label{tab: quadratic}
\end{table}

\begin{table}[t]
    \centering
    \begin{tabular}{|l|*{2}{c|}}\hline
   Global Direct & Global Dynamic & Local Dynamic \\ \hline
    \textbf{2072 s.} &
    \textbf{1309 s.}  &
    \textbf{14823 s.}   \\ \hline
    \end{tabular}
    \caption{Duration times of the methods (2000 iterations for global methods, 20000 iterations for local method) on the fully coupled quadratic model \eqref{eq: quadratic model} for $T=1$.}
\end{table}

\begin{figure}[H]
   \begin{minipage}[c]{.49\linewidth}
      \includegraphics[width=0.9\linewidth]{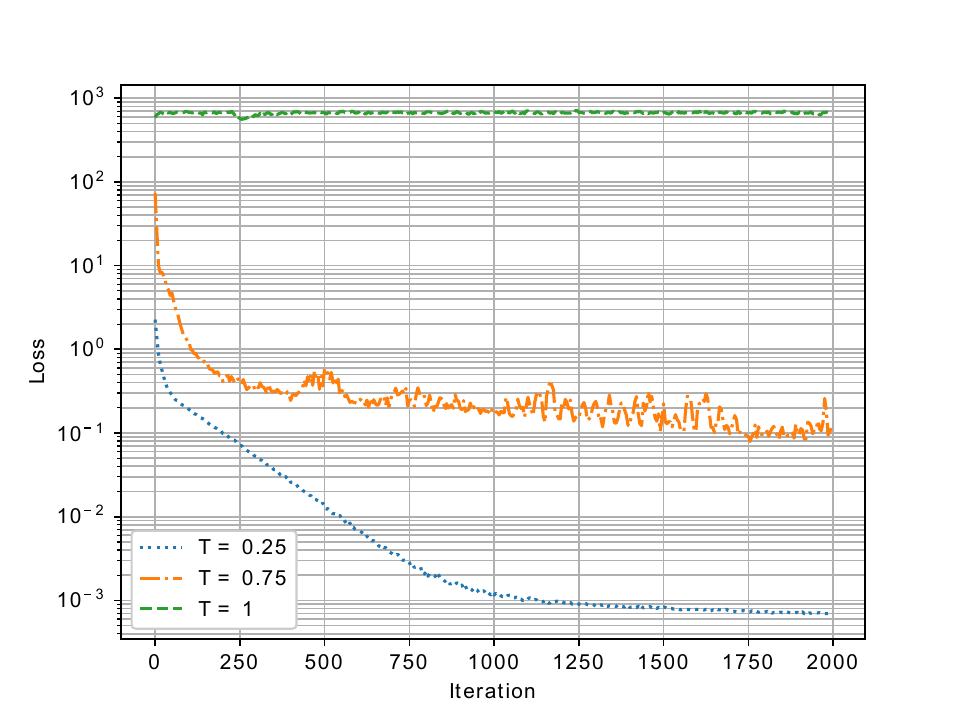}
   \end{minipage} \hfill
   \begin{minipage}[c]{.49\linewidth}
      \includegraphics[width=0.9\linewidth]{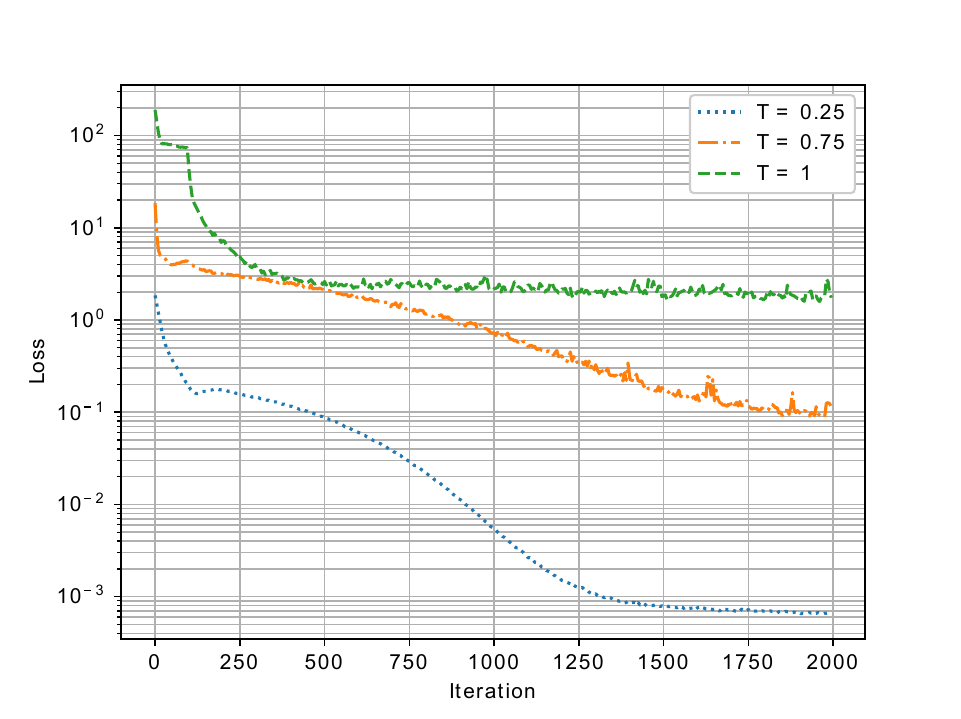}
   \end{minipage}
   \caption{Learning curves for Global Direct  (left) and Global Dynamic (right) methods on the fully coupled quadratic model \eqref{eq: quadratic model}. The loss is the $L^2$ error between $Y_T$ and the terminal condition of the backward equation.}
   \label{fig: loss 2}
\end{figure}

\begin{figure}[H]
   \centering
   \includegraphics[width=0.49\linewidth]{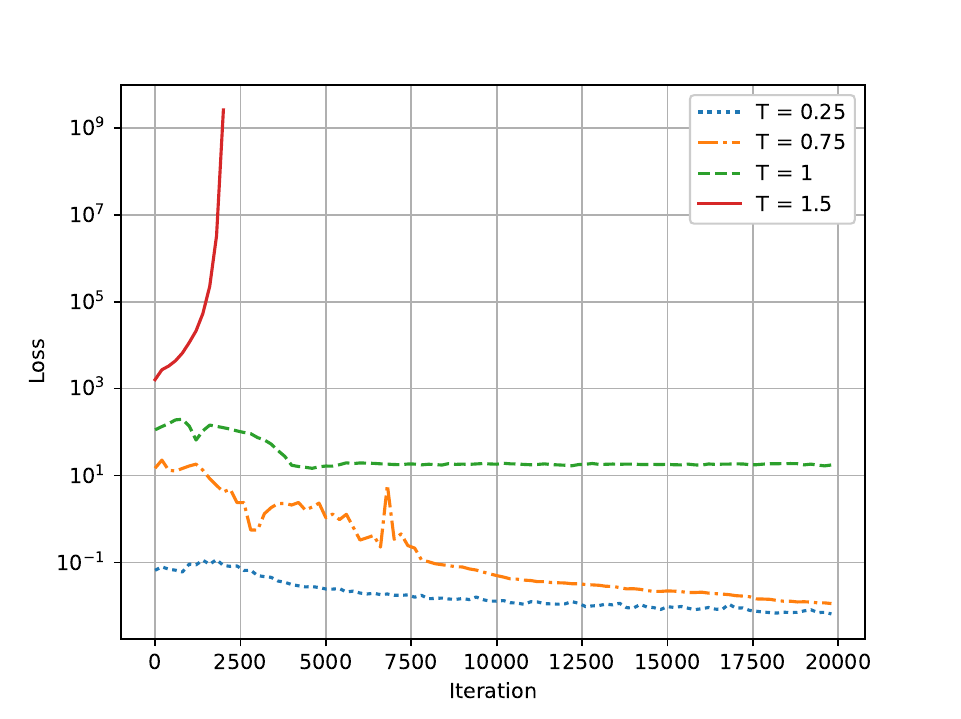}
    \caption{Learning curves for the Local Dynamic methods on the fully coupled quadratic model \eqref{eq: quadratic model}. The loss is the sum of the local $L^2$ errors between the neural network $Y$ and the Euler discretization for all time steps.}
   \label{fig: learning local quadratic}
\end{figure}

\begin{figure}[H]
    \centering 
    \includegraphics[width=0.49\linewidth]{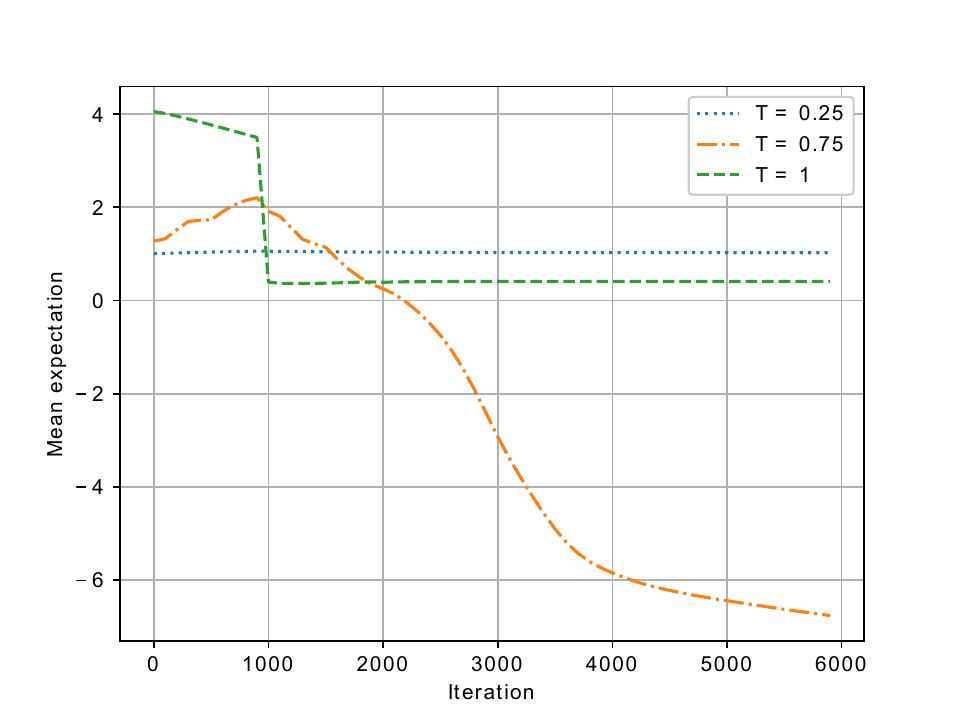}
    \caption{$\E[X_T]$ for Local Dynamic method on the fully coupled quadratic model \eqref{eq: quadratic model}.}
    \label{fig: negative exp}
\end{figure}

\begin{figure}[H]
   \begin{minipage}[c]{.49\linewidth}
      \includegraphics[width=0.9\linewidth]{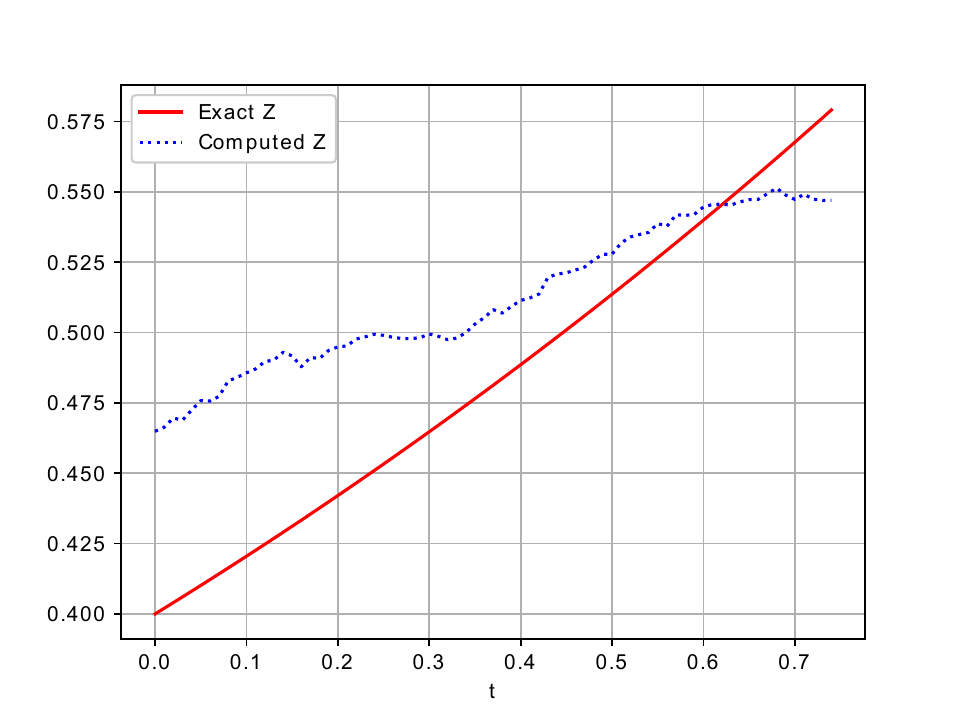}
   \end{minipage} \hfill
   \begin{minipage}[c]{.49\linewidth}
      \includegraphics[width=0.9\linewidth]{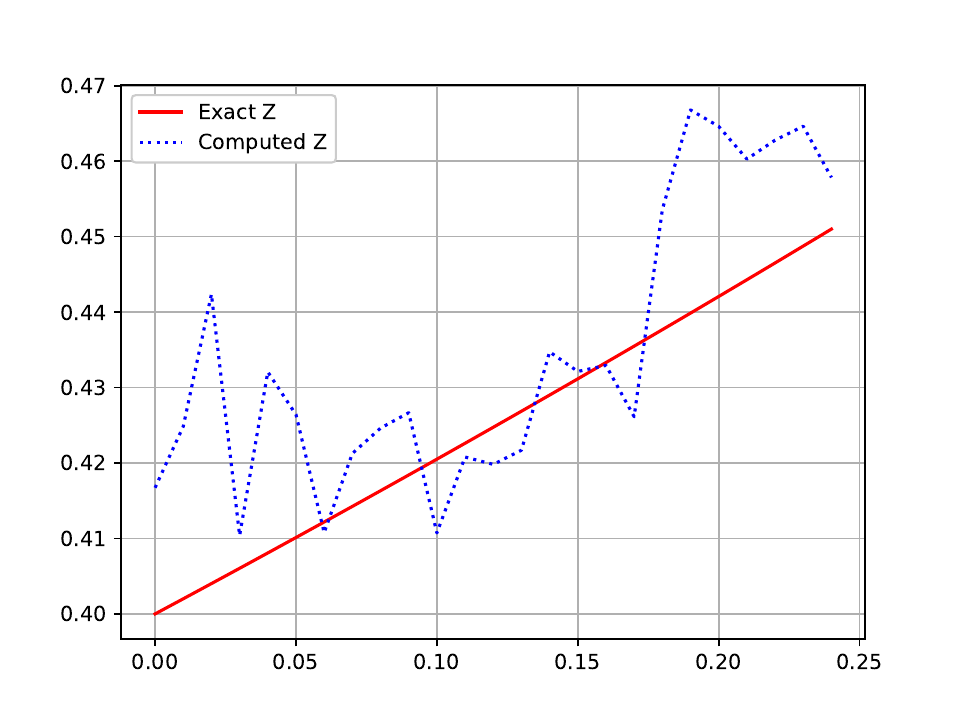}
   \end{minipage}
   \caption{First coordinate of $Z_t$ evaluated on a sample path for Global Direct  (left) ($T=0.75$)  and 
   Local Dynamic methods (right)($T=0.25$) methods after 2000 iterations (respectively 20000) on the fully coupled quadratic model \eqref{eq: quadratic model}.}
   \label{fig: Z quadratic}
\end{figure}

\begin{figure}[H]
   \begin{minipage}[c]{.49\linewidth}
      \includegraphics[width=0.9\linewidth]{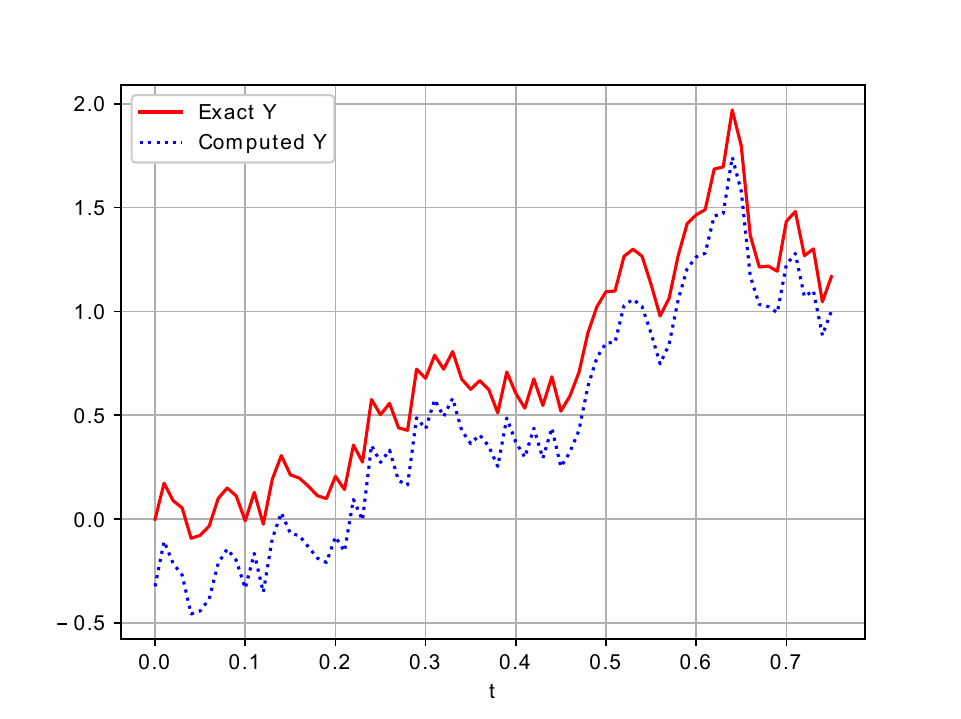}  \end{minipage} \hfill
   \begin{minipage}[c]{.49\linewidth}
      \includegraphics[width=0.9\linewidth]{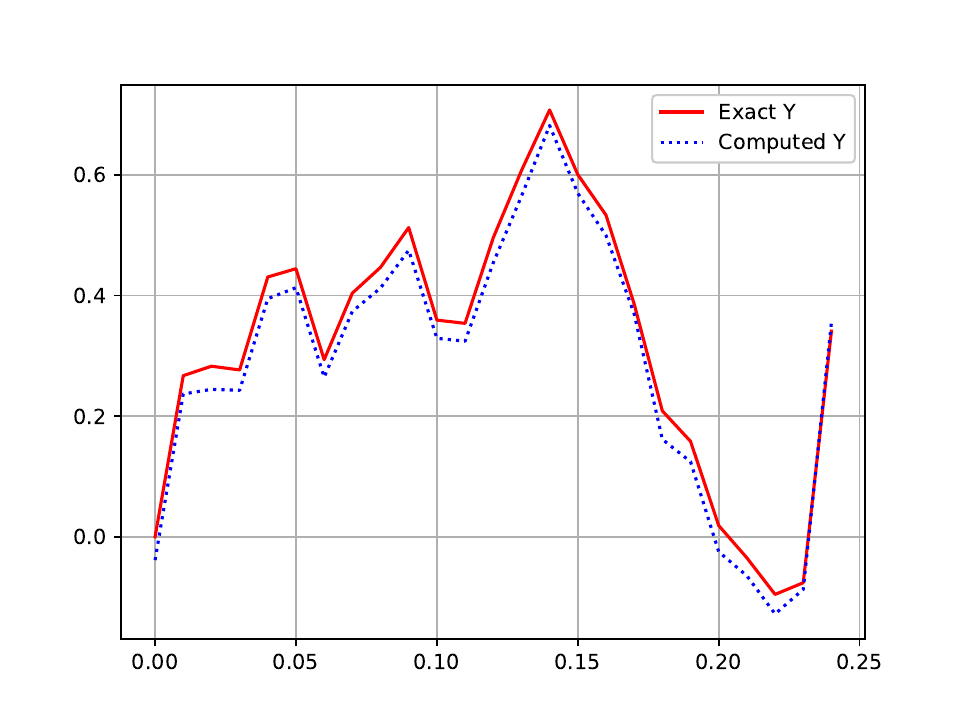}
   \end{minipage}
   \caption{First coordinate of $Y_t$ evaluated on a sample path for 
   Global Dynamic  (right) method ($T=0.75$) and Local Dynamic method (right) ($T=0.25$) after 2000 iterations (respectively 20000 for the local method) on the fully coupled quadratic model \eqref{eq: quadratic model}.}
   \label{fig: better Y}
\end{figure}

\paragraph{}
We take $a=b=c=0.1,\ \alpha = 0.5, \sigma = 0.4, \xi = 1$. We observe in Table \ref{tab: quadratic} convergence of the methods for small maturities and divergence beyond $T=1$. Note that the dynamic estimation of the expectation prevents the algorithm to explode for $T=1$, contrarily to the Global Direct method. However, it does not converge to the true solution in this case. Indeed the loss plateaus at the value 2 in Figure \ref{fig: loss 2} (right), so the terminal condition of the BSDE is not properly respected. The Global Dynamic method also produces better result than the Global Direct one (see Figure \ref{fig: better Y} and Table \ref{tab: quadratic}). We notice from Figure \ref{fig: Z quadratic} and Figure \ref{fig: better Y} that the estimated $Y,Z$ processes have the good shape but some errors are still present after convergence.

Concerning the Local Dynamic method, we see in Figure \ref{fig: negative exp} that the estimated expectations are stable around zero for a few iterations but then become negative. It may be due to the lack of a contraction for the fixed point problem. The loss explodes for $T=1.5$, as seen on the learning curve from Figure \ref{fig: learning local quadratic}. For $T=1$, we see that it stays above $10$.

\section{Conclusion}
We have shown that neural network methods can solve some moderate dimensional FBSDE of McKean-Vlasov type.
Comparing the different algorithms we find out that:
\begin{itemize}
    \item The dynamic update of the expectation is efficient in terms of computation speed (about 30\% faster than direct method) and seems to smooth the learning curve.
    \item For the mean-field games of controls example, the Pontryagin approach performs better than the Weak one for large maturities. On the contrary, the Weak approach is the best for small maturities.
    \item For the fully coupled linear model we observe no convergence problem whereas for the fully coupled quadratic one we can solve only the equation on a small time horizon. However the Local Dynamic method is not very accurate for larger maturities.
    \item The Local Dynamic method faces more difficulties for quadratic problems than the global methods do. It also requires more iterations, hence more time, to converge.
    \item The expectation law estimation method does not work well and requires to empirically choose a proper penalization parameter, which is troublesome.
    \item The methods can be used in dimension 10, thus applied to more realistic problems than usually. For instance, in the price impact model, the number of dimensions corresponds to the number of assets involved in the trading. Thus, developing methods able to deal with problems in high dimensions can help us to handle large portfolios.
\end{itemize}
We recommend the use of the global method combined with dynamic moment estimation which offers the best accuracy and training speed among all the tested methods. For linear quadratic mean-field games of controls, it appears to be better to use the Weak approach for small maturities and the Pontryagin method for larger time horizons. The use of a Local Dynamic method is possible but requires too many iterations to converge hence it is not competitive in terms of computation time.

\printbibliography

\end{document}